\documentclass[11pt]{article}
\usepackage{amsmath}
\usepackage{pdfsync}
\usepackage{hyperref}
\usepackage{ifpdf}
\usepackage{amsfonts}
\usepackage{amssymb}
\usepackage{graphicx}
\usepackage[monochrome]{color}
\newtheorem{thm}{Theorem}[section]
\newtheorem{cor}[thm]{Corollary}
\newtheorem{lem}[thm]{Lemma}
\newtheorem{prop}[thm]{Proposition}
\newtheorem{defn}[thm]{Definition}

\numberwithin{equation}{section}

\newcommand{\dx}{\,{\rm d}x}
\newcommand{\dy}{\,{\rm d}y}

\newcommand{\dt}{\,{\rm d}t}
\newcommand{\rd}{{\rm d}}

\def\LL{\mathrm{L}} 
\def\supp{\mathrm{supp}} 
\newcommand{\RR}{\mathbb{R}}

\newcommand{\ren}{\RR^d}
\def\ee{\mathrm{e}} 
\def\qed{\,\unskip\kern 6pt \penalty 500
\raise -2pt\hbox{\vrule \vbox to8pt{\hrule width 6pt
\vfill\hrule}\vrule}\par}
\definecolor{darkblue}{rgb}{0.05, .05, .65}
\definecolor{darkgreen}{rgb}{0.1, .65, .1}
\definecolor{darkred}{rgb}{0.8,0,0}

\begin{document}
\title{\textbf{ Quantitative Local and Global  A Priori Estimates \\ for Fractional Nonlinear Diffusion Equations}\\[7mm]}

\author{\Large Matteo Bonforte$^{\,a,\,b}$ 
~and~ Juan Luis V\'azquez$^{\,a,\,c}$\\} 
\date{} 

\maketitle

\

\begin{abstract}

We establish  quantitative estimates for  solutions
$u(t,x)$ to the fractional nonlinear diffusion equation, $\partial_t u
+(-\Delta)^s (u^m)=0$ \ in the whole range of exponents $m>0$, $0<s<1$. The equation is posed in the whole space $x\in\RR^d$. We first obtain  weighted global integral estimates that allow to establish existence of solutions for classes of large data. In the core of the paper we obtain quantitative pointwise lower estimates of the positivity of the solutions, depending only on the norm of the initial data in a certain ball.  The estimates take a different form in three exponent ranges: slow diffusion, good range of fast diffusion, and very fast diffusion. Finally, we show existence and uniqueness of initial traces.
\end{abstract}

\vskip 1cm

\vspace{1cm}

\begin{itemize}
\item[(a)] Departamento de Matem\'{a}ticas, Universidad
Aut\'{o}noma de Madrid,\\ Campus de Cantoblanco, 28049 Madrid, Spain
\item[(b)] e-mail address:~\texttt{matteo.bonforte@uam.es }\\ web-page:~\texttt{http://www.uam.es/matteo.bonforte}
\item[(c)] e-mail address:~\texttt{juanluis.vazquez@uam.es }\\ web-page:~\texttt{http://www.uam.es/juanluis.vazquez}
\end{itemize}

\newpage

\tableofcontents

\newpage

\section{Introduction}

We consider the class of nonnegative weak solutions of the fractional diffusion equation
\begin{equation}\label{FDE.eq}
\partial_t u +(-\Delta)^s (u^m)=0\quad \; \mbox{in }(0,T)\times\RR^d\,,
\end{equation}
where $m>0$, $0<s<1$, $d\ge 1$, and $T>0$. The precise definition of the fractional Laplacian is given in Appendix \ref{ssec.app1}.   We are mainly interested in values $m\ne1$, since the linear case is rather well known. For $s=1$ we recover the classical porous medium/fast diffusion equation (that will be shortened as PME/FDE respectively), whose theory is well-known, cf. \cite{VazBook}. We will call the case $s=1$ the standard diffusion case. \normalcolor We also assume that we are given initial data
\begin{equation}\label{FDE.id}
u(0,x)=u_0(x)\,,
\end{equation}
where in principle $u_0\in \LL^1(\RR^d)$ and $u_0\ge 0$. However, larger classes of initial data are sometimes considered, like changing sign solutions or non integrable data. The equation has recently attracted some attention in mathematical analysis. Such an interest  has been motivated by its appearance as a model for anomalous diffusion in different applied contexts. For the reader's convenience we have listed in Appendix II (Section \ref{App.Motiv}) the most relevant sources to the applications that we know of. \normalcolor

We refer to \cite{DPQRV1,DPQRV2} for the basic theory of existence and uniqueness of weak solutions for the Cauchy problem  \eqref{FDE.eq}$-$\eqref{FDE.id}.   These papers also  describe the main results on $L^q$ boundedness and $C^\alpha$ regularity, they show that nonnegative solutions are indeed positive everywhere, as well as some other basic properties of the nonlinear semigroup generated by the problem. Recently, the existence and properties of Barenblatt solutions for the Cauchy Problem was established in \cite{Vaz2012}.  Related literature is also mentioned in these papers. \normalcolor

The main purpose of the present paper is obtaining quantitative a  priori estimates of a local type for the solutions of the problem. Such estimates were obtained for the standard PME by Aronson-Caffarelli \cite{ArCaff} and by the authors for the standard FDE \cite{BV, BV3, BV-ADV}\,. This is not always possible for the present model due to the nonlocal character of the diffusion operator, but then global estimates occur in weighted spaces. The use of suitable weight functions allows to prove
crucial $\LL^1$-weighted estimates that enter substantially into the derivation of the main results. The results take different forms according to the value of the exponent $m$, a fact that is to be expected since it happens for standard diffusion. The list of bounds is as follows: weighted $\LL^1$ estimates, half-Harnack parabolic estimates (i.\,e., quantitative pointwise lower estimates), and tail estimates (i.\,e., asymptotic spatial behaviour). As a first consequence of these estimates, an existence result for very weak solutions with non-integrable data in some weighted $\LL^1_{\varphi}$ space is obtained. In particular, bounded initial data or data with slow growth at infinity are allowed.

We remind that Harnack inequalities are a standard tool used to develop further theory, see in that respect the work of Di Benedetto for quasilinear equations \cite{DBbook, DGVbook}. Let us quote four examples of application of the line of results of this paper that are already available: the existence and uniqueness of initial traces, which we do in Section 7; understanding the asymptotic behaviour of the fractional KPP equation \cite{SV2013}; dealing with nonlocal symmetrization problems \cite{VV2013}; or uniqueness issues for the fractional PME with variable density \cite{PT2}.\normalcolor

\medskip

\textbf{Outline of the paper and main results. }
  First, some preliminary information on case divisions. The case $m>1$ is called the {\sl (fractional) porous medium case} :  contrary to the standard porous medium equation, it does not have the property of finite propagation, an important difference established in \cite{DPQRV1,DPQRV2}\,. The range of exponents $m\in (0,1)$ is called the {\sl (fractional) fast diffusion equation}, and it has special properties when $(d-2s)/d=: m_c<m<1$\,, which we call the {\sl good fast diffusion range}. When $0<m<m_c$ it is known that some solutions extinguish in finite time, which is a clear manifestation of the change of character of the equation, since solutions of the Cauchy problem exist globally in time and are positive everywhere in $Q=(0,+\infty)\times\RR^d$ if $m\ge m_c$\,.\normalcolor

In Section \ref{sect.2}, we derive integral bounds in form of weighted $\LL^1$ estimates, valid for nonnegative solutions of the Cauchy problem in the whole fast diffusion range $0<m<1$. Actually, they are valid for the difference of two ordered solutions, the precise statement is given in Theorem \ref{prop.HP.s}. Contrary to the purely local $\LL^1$ estimates known in the standard fast diffusion case, cf. \cite{HP}, the estimates for $s<1$ are valid in weighted $\LL^1$-spaces and the weight must  decay at infinity with a certain decay rate, not too fast, not too slow. This is again a manifestation of the nonlocal properties of the fractional Laplacian. The estimates will be important as a priori bounds for solutions, or families of solutions, through the rest of the paper.

In Section \ref{sect.exist.large} we use the estimates of  Section \ref{sect.2} to construct solutions for initial data that belong to weighted $\LL^1$-spaces, in particular for data $0\le u_0\in \LL^1_{\rm loc}(\RR^d)$ such that $u_0(x)$ grows less than $O(|x|^{2s/m})$ as $|x|\to\infty$, in particular for all bounded data. These solutions can be uniquely identified as minimal solutions in a precise sense and satisfy many of the properties of the known class of bounded and integrable weak solutions.

Section \ref{sect.3} studies the actual positivity of nonnegative solutions via quantitative lower estimates for the good fast diffusion equation. Precise local lower bounds are contained in Theorem \ref{thm.lower}. The behaviour as $|x|\to\infty$ (so-called tail behaviour) is studied in Section \ref{ssec.gfd.tail}, and global spatial lower bounds are derived as a consequence in Section \ref{ssec.global.GFDE}. The merit of the estimates is that they are quantitative and most of the exponents are sharp.  The lower estimates of this section can be adapted for the exponent $m=m_c$ separating both fast diffusion subranges, but only when $u_0\in \LL^p_{loc}$ for some $p>1$. However, we refrain from doing this particular case in the present paper since the proof uses some other techniques that would lengthen the text.

The very fast diffusion range $0<m<m_c$ is studied in Section \ref{sect.4}. The weighted $\LL^1$ estimates of Theorem \ref{prop.HP.s} continue to hold, but this does not allow to obtain the same type of  quantitative lower bounds since the technique used in the good fast diffusion range does not work anymore. There are two problems: on the one hand the $L^1$--$L^\infty$ smoothing effect does not hold for general $\LL^1$ initial data, on the other hand the presence of the extinction phenomenon makes things more complicated, and the extinction time enters directly the estimates of Theorem \ref{thm.lower.subcrit.}. These difficulties have already appeared in the standard FDE, $s=1$, and were treated in our paper \cite{BV-ADV}. However the technique used in that paper does not extend to $0<s<1$ and we present here a technique that is based on the careful use of weight factors, and in the limit $s=1$ gives a simpler proof of the result of \cite{BV-ADV}. We also study the problem of characterizing the finite time extinction in terms of the initial data; thus, we determine a class of initial data that produces solutions that extinguish in finite time, see Proposition \ref{prop.ext}, as well as a roughly complementary class of initial data for which the solution exists and is positive globally in time, see Corollary \ref{cor.not.ext.}\,.

Section \ref{sect.PME} is devoted to study similar questions for the porous medium case. Theorem \ref{thm.lower.pme} establishes local lower bounds of the Aronson-Caffarelli type for all $0<s\le 1$. The question of optimal decay as $|x|\to\infty$ is an open problem; for selfsimilar solutions it is solved in \cite{Vaz2012}.

In Section \ref{sect.traces} we address a different question that complements our previous results, i.\,e., the question of existence and uniqueness of an initial trace for nonnegative weak solutions defined in a strip $Q_T=(0,T]\times\RR^d$. The main results are stated in Theorems \ref{thm.init.trace.m<1} and \ref{thm.init.trace.m>1}. This result can be combined in the reverse direction with the existence of solutions with initial data a nonnegative Radon measure, Theorem 4.1 of \cite{Vaz2012}.

In  Appendix I we collect the definitions of weak, very weak and strong solutions, together with a number of technical results. As already mentioned above, Appendix II discusses applications.

\medskip

\noindent {\sc About the linear equation.} Our estimates have counterparts for the linear fractional heat equation, case $m=1$, that are worth commenting. Thus, the lower bound of Section \ref{sect.PME} for $m>1$ passes to the limit $m\downarrow 1$, and this coincides with the limit $m\uparrow1$ of a part of the estimate for $m<1$ obtained in Section \ref{sect.3} for $m<1$. See Proposition \ref{thm.lower.m1}.
Finally, we prove the existence and uniqueness of an initial trace also when $m=1$, cf. Theorem \ref{thm.init.trace.m=1}. This is an interesting result that is not present in the literature to our knowledge and complements the uniqueness results of \cite{BPSV2013}. \normalcolor

\medskip

\noindent{\sc Notations.}  Throughout the paper, we fix $m_c=(d-2s)/d$, $m_1=d/(d+2s)$, $p_c=d(1-m)/2$, and $\vartheta:=1/[2s-d(1-m)]$\,, which is positive if $m>m_c$\,. We will call $s$-Laplacian of $f$ the function $-(-\Delta)^{s}f$. This is consistent with the use in the standard case $s=1$.\normalcolor

\section{Weighted $\LL^1$ estimates in the fast Diffusion range}\label{sect.2}

We  will derive weighted $\LL^1$ estimates which also hold for the standard FDE (i.e., the limit case $s=1$). When $s<1$ the equation is nonlocal, therefore we cannot expect purely local estimates to hold. Indeed we will obtain estimates in weighted spaces if the weight satisfies certain decay conditions at infinity. We present first a technical lemma which will be used several times in the rest of the paper.

\begin{lem}\label{Lem.phi}
Let $\varphi\in C^2(\RR^d)$ and positive real function that is radially symmetric and decreasing in $|x|\ge 1$. Assume also that $\varphi(x)\le |x|^{-\alpha}$ and that $|D^2\varphi(x)| \le c_0 |x|^{-\alpha-2}$\,, for some positive constant $\alpha$ and for $|x|$ large enough. Then,  for all $|x|\ge |x_0|>>1$ we have
\begin{equation}\label{Delta.s.phi}
|(-\Delta)^s\varphi(x)|\le
\left\{\begin{array}{lll}
\dfrac{c_1}{|x|^{\alpha+2s}}\,,   & \mbox{if $\alpha<d$}\,,\\[5mm]
\dfrac{c_2\log|x|}{|x|^{d+2s}}\,,   & \mbox{if $\alpha=d$}\,,\\[5mm]
\dfrac{c_3}{|x|^{d+2s}}\,,   & \mbox{if $\alpha>d$}\,,\\[5mm]
\end{array}\right.
\end{equation}
with positive constants $c_1,c_2,c_3>0$ that depend only on $\alpha,s,d$ and $\|\varphi\|_{C^2(\RR^d)}$. For  $\alpha>d$ the reverse estimate holds  from below  if $\varphi\ge0$: \ $|(-\Delta)^s\varphi(x)|\ge c_4 |x|^{-(d+2s)}$  for all $|x|\ge |x_0|>>1$\,.
\end{lem}

The proof is easy but technical, and is given in Appendix \ref{sec.A1} for the reader's convenience.
We point out that the large-decay case $\alpha>d$ is what makes the estimate in the fractional Laplacian case very different from the usual Laplacian case. In particular, the $s$-Laplacian of a nonnegative smooth function with compact support is strictly positive outside of the support and has a certain decay at infinity, indeed the minimal decay $|x|^{-(d+2s)}$ is obtained for the $(-\Delta)^{s}\varphi$ when $\varphi\ge 0$ is compactly supported, cf. \cite{DPQRV2}.
A suitable particular choice is the function $\varphi$ defined for $\alpha>0$ as $\varphi(x)=1$ for $|x|\le 1$ and
\begin{equation}\label{phi}
\varphi(x)=
\dfrac{1}{\left(1+(|x|^2-1)^4\right)^{\alpha/8}}\,, \qquad\mbox{if } |x|\ge 1\,.
\end{equation}

We are now ready to present the weighted estimates.

\begin{thm}[Weighted $\LL^1$ estimates]\label{prop.HP.s}
Let $u\ge v$ be two ordered solutions to the  equation \eqref{FDE.eq}, with $0<m<1$. Let $\varphi_R(x)=\varphi(x/R)$ where $R>0$ and $\varphi$ is as in the previous lemma with $0\le \varphi(x)\le |x|^{-\alpha}$ for $|x|>>1$ and
$$
d-\frac{2s}{1-m}<\alpha< d+\frac{2s}{m}\,.
$$
Then, for all $0\le \tau,t <\infty$ we have
\begin{equation}\label{HP.s}
\left(\int_{\RR^d}\big(u(t,x)- v(t,x)\big)\varphi_R(x)\dx\right)^{1-m}\le
\left(\int_{\RR^d}\big(u(\tau,x)- v(\tau,x)\big)\varphi_R(x)\dx\right)^{1-m}
+ \frac{C_1 \,|t-\tau|}{R^{2s-d(1-m)}}
\end{equation}
with $C_1>0$ that depends only on $\alpha,m,d$\,.
\end{thm}

It is remarkable that the  estimate holds for (very) weak solutions, maybe changing sign.  Also, it is worth pointing out that the estimate holds both for $\tau<t$ and for $\tau>t$. In the limit $s\to 1$ we recover the well known $\LL^1$ local estimates for the standard FDE.

\noindent {\sl Proof.~} \noindent$\bullet~$\textsc{Step 1. }\textit{A differential inequality for the weighted $\LL^1$-norm. }If $\psi$ is a smooth and sufficiently decaying function we have
\[
\begin{split}
\left|\frac{\rd}{\dt}\int_{\RR^d}\big(u(t,x)-v(t,x)\big)\psi(x)\dx\right|
&=\left|\int_{\RR^d}\left((-\Delta)^s u^m-(-\Delta)^s v^m\right)\psi\dx\right|\\
&=_{(a)}\left|\int_{\RR^d}\left(u^m- v^m\right)(-\Delta)^s\psi\dx\right|\\
&\le_{(b)} 2^{1-m}\int_{\RR^d}(u- v)^m\left|(-\Delta)^s\psi\right|\dx\\
&\le_{(c)} 2 \left(\int_{\RR^d}(u- v)\psi\dx\right)^m\,
        \left(\int_{\RR^d}\frac{\left|(-\Delta)^s\psi\right|^{\frac{1}{1-m}}}{\psi^\frac{m}{1-m}}\dx\right)^{1-m}\,.
\end{split}
\]
Notice that in $(a)$ we have used the fact that $(-\Delta)^s$ is a symmetric operator, while in $(b)$ we have used that $\left(u^m- v^m\right)\le 2^{1-m}(u-v)^m$, where $u^m=|u|^{m-1}u$ as mentioned. In $(c)$ we have used H\"older inequality with conjugate exponents $1/m>1$ and $1/(1-m)$. If the last integral factor is bounded, then we get
\[
\left|\frac{\rd}{\dt}\int_{\RR^d}\big(u(t,x)-v(t,x)\big)\psi(x)\dx\right|\le C_\psi^{1-m}\left(\int_{\RR^d}\big(u(t,x)- v(t,x)\big)\psi(x)\dx\right)^m
\]
Integrating the above differential inequality on $(\tau,t)$ with $\tau,t\ge 0$ we obtain:
\[
\left(\int_{\RR^d}\big(u(t,x)- v(t,x)\big)\psi(x)\dx\right)^{1-m}-\left(\int_{\RR^d}\big(u(\tau,x)- v(\tau,x)\big)\psi(x)\dx\right)^{1-m}
\le (1-m)C_\psi^{1-m}\,|t-s|
\]
which is \eqref{HP.s} once we estimate the constant $C_\psi$, for a convenient choice of test function.

\noindent$\bullet~$\textsc{Step 2. }\textit{Estimating the constant $C_\psi$. }Choose $\psi(x)=\varphi_R(x):=\varphi(x/R)=\varphi(y)$\,, with $\varphi$ as in Lemma \ref{Lem.phi} and $y=x/R$\,, so that $(-\Delta)^s\psi(x)=(-\Delta)^s\varphi_R(x)=R^{-2s}(-\Delta)^s\varphi(y)$
\[\begin{split}
C_\psi
&=\int_{\RR^d}\frac{\left|(-\Delta)^s\varphi_R(x)\right|^{\frac{1}{1-m}}}{\varphi_R(x)^\frac{m}{1-m}}\dx
 =R^{d-\frac{2s}{1-m}}\int_{\RR^d}\frac{\left|(-\Delta)^s\varphi(y)\right|^{\frac{1}{1-m}}}{\varphi(y)^\frac{m}{1-m}}\dy\\
&=R^{d-\frac{2s}{1-m}}\left[\int_{B_{2}}\frac{\left|(-\Delta)^s\varphi(y)\right|^{\frac{1}{1-m}}}{\varphi(y)^\frac{m}{1-m}}\dy
    +\int_{B_{2}^{\,c}}\frac{\left|(-\Delta)^s\varphi(y)\right|^{\frac{1}{1-m}}}{\varphi(y)^\frac{m}{1-m}}\dy\right]
    = k_1 R^{d-\frac{2s}{1-m}}\,,\\
\end{split}\]
where it is easy to check that the first integral is bounded, since $\varphi\ge k_2>0$ on $B_{2}$\,, and when $|y|>|x_0|$ with $|x_0|>>1$ we know by estimates \eqref{Delta.s.phi} that
\begin{equation}\label{c.psi}
\frac{\left|(-\Delta)^s\varphi(y)\right|^{\frac{1}{1-m}}}{\varphi(y)^\frac{m}{1-m}}\le
\left\{\begin{array}{lll}
\dfrac{k_3}{|y|^{\alpha+\frac{2s}{1-m}}}\,,   & \mbox{if $\alpha<d$}\,,\\[5mm]
\dfrac{k_4\log|y|}{|y|^{d+\frac{2s}{1-m}}}\,,   & \mbox{if $\alpha=d$}\,,\\[5mm]
\dfrac{k_5}{|y|^{\frac{d+2s-\alpha m}{1-m}}}\,,   & \mbox{if $\alpha>d$}\,,\\[5mm]
\end{array}\right.
\end{equation}
therefore $k_1$ is finite whenever $d-\frac{2s}{1-m}<\alpha< d+\frac{2s}{m}$\,. Note that all the constants $k_i$ depend only on $\alpha, m, d$\,.\qed

\noindent\textbf{Remark. }The estimate implies the conservation of mass when $(d-2s)/d=m_c<m<1$, by letting $R\to \infty$.
On the other hand,  when $0<m<m_c$ solutions corresponding to $u_0\in\LL^1(\RR^d)\cap\LL^p(\RR^d)$ with $p\ge d(1-m)/2s$\,, extinguish in finite time $T>0$\,, (see e.g. \cite{DPQRV2}); the above estimates provide a lower bound for the extinction time in such a case, just by letting $\tau=T$ and $t=0$ in the above estimates:
\begin{equation}\label{HP.s.T}
\frac{1}{C_1\,R^{d(1-m)-2s}}\left(\int_{\RR^d}u_0\,\varphi_R\dx\right)^{1-m}\le T
\end{equation}
Moreover, if the initial datum $u_0$ is such that the limit as $R\to+\infty$ of the right-hand side diverges to $+\infty$, then the corresponding solution $u(t,x)$ exists (and is positive) globally in time, as explained in Corollary \ref{cor.not.ext.}\,.

\section{Existence of solutions in weighted $\LL^1$-spaces}\label{sect.exist.large}

\begin{thm}\label{exist.large}
Let $0<m<1$ and let $u_0\in\LL^1(\RR^d, \varphi\dx)$, where $\varphi$ is as in Theorem $\ref{prop.HP.s}$ with decay at infinity $|x|^{-\alpha}$, $d-[2s/(1-m)]<\alpha<d+(2s/m)$. Then there exists a  very weak solution  $u(t,\cdot)\in\LL^1(\RR^d, \varphi\dx)$ to equation \eqref{FDE.eq} on $[0,T]\times \RR^d$, in the sense that
\[
\int_0^T\int_{\RR^d}u(t,x)\psi_t(t,x)\dx\dt
=\int_0^T\int_{\RR^d}u^m(t,x)(-\Delta)^s\psi(t,x)\dx\dt\,,\qquad\mbox{for all $\psi\in C_c^\infty([0,T]\times\RR^d)\,.$}
\]
This solution is continuous in the weighted space, $u\in C([0,T]:\LL^1(\RR^d, \varphi\dx))$\,.
\end{thm}

\noindent {\sl Proof.~}Let $\varphi=\varphi_R$ be as in Theorem \ref{prop.HP.s} with the decay at infinity $|x|^{-\alpha}$. Let $0\le u_{0,n}\in \LL^1(\RR^d)\cap \LL^\infty(\RR^d)$ be a non-decreasing sequence of initial data $u_{0,n-1}\le u_{0,n}$, converging monotonically to $u_0\in \LL^1(\RR^d, \varphi\dx)$\,, i.\,e., such that $\int_{\RR^d}(u_0- u_{n,0})\varphi\dx \to 0$ as $n\to \infty$. Consider the unique solutions $u_n(t,x)$ of equation \eqref{FDE.eq} with initial data $u_{0,n}$. By the comparison results of \cite{DPQRV2} we know that they form a monotone sequence. The weighted estimates \eqref{HP.s} show that the sequence is bounded in $\LL^1(\RR^d, \varphi\dx)$ uniformly in $t\in[0,T]$\,. By the monotone convergence theorem in $\LL^1(\RR^d, \varphi\dx)$, we know that the solutions $u_n(t,x)$ converge monotonically as $n\to \infty$ to a function $u(t,x)\in \LL^\infty ((0,T): \LL^1(\RR^d, \varphi\dx))$. Indeed, the weighted estimates \eqref{HP.s} show that when $u_0\in \LL^1(\RR^d, \varphi\dx)$ then
\begin{equation}\label{HP.s.3}\begin{split}
\left(\int_{\RR^d}u(t,x)\varphi(x)\dx\right)^{1-m}
&=\lim_{n\to\infty}\left(\int_{\RR^d}u_n(t,x)\varphi(x)\dx\right)^{1-m}\\
&\le \lim_{n\to\infty} \left(\int_{\RR^d}\big(u_n(0,x)\big)\varphi(x)\dx\right)^{1-m} + C_1 R^{d(1-m)-2s}\,t\\
&=\left(\int_{\RR^d}u_0(x)\varphi(x)\dx\right)^{1-m}
+ C_1 R^{d(1-m)-2s}\,t
\end{split}
\end{equation}
At this point we need to show that the function $u(t,x)$ constructed as above is a very weak solution to equation \eqref{FDE.eq} on $[0,T]\times \RR^d$\,, more precisely we have to show that for all $\psi\in C_c^\infty([0,T]\times\RR^d)$ we have
\begin{equation}\label{vw.sol}
\int_0^T\int_{\RR^d}u(t,x)\psi_t(t,x)\dx\dt
=\int_0^T\int_{\RR^d}u^m(t,x)(-\Delta)^s\psi(t,x)\dx\dt\,.
\end{equation}
By the results of \cite{DPQRV2} we know that each $u_n$ is a bounded strong solutions, since the initial data $u_0\in \LL^1(\RR^d)\cap\LL^\infty(\RR^d)$\,, therefore for all $\psi\in C_c^\infty([0,T]\times\RR^d)$ we have
\begin{equation}\label{sol.approx}
\int_0^T\int_{\RR^d}u_n(t,x)\psi_t(t,x)\dx\dt
=\int_0^T\int_{\RR^d}u^m_n(t,x)(-\Delta)^s\psi(t,x)\dx\dt\,.
\end{equation}
Now, for any $\psi\in C_c^\infty([0,T]\times\RR^d)$ we easily have that
\[
\lim_{n\to\infty}\int_0^T\int_{\RR^d}u_n(t,x)\psi_t(t,x)\dx=\int_0^T\int_{\RR^d}u(t,x)\psi_t(t,x)\dx
\]
since $\psi$ is compactly supported and we already know that $u_n(t,x)\to u(t,x)$ in $L^1_{\rm loc}$.
On the other hand, for any $\psi\in C_c^\infty([0,T]\times\RR^d)$ we have that
\[
\lim_{n\to\infty}\int_0^T\int_{\RR^d}u^m_n(t,x)(-\Delta)^s\psi(t,x)\dx\dt=\int_0^T\int_{\RR^d}u^m(t,x)(-\Delta)^s\psi(t,x)\dx\dt
\]
since $u_n\le u$ and
\[\begin{split}
0\le & \int_0^T\int_{\RR^d}(u^m(t,x)-u_n^m(t,x))(-\Delta)^s\psi(t,x)\dx\dt\\
&\le \int_0^T\int_{\RR^d}|u(t,x)-u_n(t,x)|^m\,\varphi^m(x)\frac{|(-\Delta)^s\psi(t,x)|}{\varphi^m(x)}\dx\dt\\
&\le \int_0^T\left(\int_{\RR^d}|u(t,x)-u_n(t,x)|\,\varphi(x)\dx\right)^m
    \left(\int_{\RR^d}\left|\frac{|(-\Delta)^s\psi(t,x)|}{\varphi(x)^m}\right|^{\frac{1}{1-m}}\dx\dt\right)^{1-m}\\
&\le C \int_0^T\int_{\RR^d}(u(t,x)-u_n(t,x))\varphi\dx\dt\to 0
\end{split}\]
where we have used H\"older inequality with conjugate exponents $1/m$ and $1/(1-m)$, and we notice that
\[
\left(\int_{\RR^d}\left|\frac{|(-\Delta)^s\psi(t,x)|}{\varphi(x)^m}\right|^{\frac{1}{1-m}}\dx\dt\right)^{1-m}\le C
\]
since $\psi$ is compactly supported, therefore by Lemma \ref{Lem.phi} we know that $\left|(-\Delta)^s\psi(t,x)\right|\le c_3|x|^{-(d+2s)}$, and the quotient
\[
\left|\frac{|(-\Delta)^s\psi(t,x)|}{\varphi(x)^m}\right|^{\frac{1}{1-m}}\le \frac{c_3}{|x|^{\frac{d+2s-m \alpha}{1-m}}}
\]
is integrable when $\frac{d+2s-m \alpha}{1-m}>d$ that is when $\alpha<d+(2s/m)$. In the last step we already know that $\int_{\RR^d}(u(t,x)-u_n(t,x))\varphi\dx\to 0$ when $\varphi$ is as above, i.e. as in Theorem \ref{prop.HP.s}. Therefore we can let $n\to\infty$ in \eqref{sol.approx} and obtain \eqref{vw.sol}\,.

For the solutions constructed above, the weighted estimates \eqref{HP.s} show that when $0\le u_0\in \LL^1(\RR^d, \varphi\dx)$ imply
\begin{equation}\label{HP.s.44}
\left|\int_{\RR^d}u(t,x)\varphi_R(x)\dx- \int_{\RR^d}u(\tau,x)\varphi_R(x)\dx\right|\le
2^{\frac{1}{1-m}}C_1 R^{d-\frac{2s}{1-m}}\,|t-\tau|^{\frac{1}{1-m}}
\end{equation}
which gives the continuity in $\LL^1(\RR^d, \varphi\dx)$\,.
Therefore, the initial trace of this solution is given by $u_0\in\LL^1(\RR^d, \varphi\dx)\,.$ \qed

\medskip

\noindent\textbf{Remark. }The solutions constructed above only need to be integrable with respect to the weight $\varphi$, which has a tail of order less than $d+2s/m$. Therefore, we have proved existence of solutions corresponding to initial data $u_0$ that can grow at infinity as $|x|^{(2s/m)-\varepsilon}$ for any $\varepsilon >0$\,. Note that for the linear case $m=1$ this exponent is optimal in view of the representation of solutions in terms of the fundamental solution, but this does not seem to be the case for $m<1$.

\medskip

\begin{thm}[Uniqueness]
The solution constructed in Theorem $\ref{exist.large}$ by approximation from below is unique. We call it the minimal solution. In this class of solutions the standard comparison result holds, and also the estimates of Theorem $\ref{prop.HP.s}$\,.
\end{thm}

\noindent {\sl Proof.~} We keep the notations of the proof of Theorem \ref{exist.large}.   Assume that there exist another sequence $0\le v_{0,k}\in \LL^1(\RR^d)$ which is monotonically non-decreasing and converges monotonically to $u_0\in \LL^1(\RR^d, \varphi\dx)$\,. By the same considerations as in the proof of Theorem \ref{exist.large}, we can show that there exists a solution $v(t,x)\in C([0,T]:\LL^1(\RR^d, \varphi\dx))$. We want to show that $u=v$, where $u$ is the solution constructed in the same way from the sequence $u_{0,n}$. We will prove equality by proving that $v\le u$ and then that $u\le v$. To prove that $v\le u$ we use the estimates
\begin{equation}\label{final.001}
\int_{\RR^d}\big[v_k(t,x)-u_n(t,x)\big]_+\dx\le \int_{\RR^d}\big[v_k(0,x)-u_n(0,x)\big]_+\dx
\end{equation}
which hold for any $u_n(t,\cdot), v_k(t,\cdot)\in \LL^1(\RR^d)$, see Theorem 6.2 of \cite{DPQRV2} for a proof. Letting $n\to \infty$ we get that
\[
\lim_{n\to\infty}\int_{\RR^d}\big[v_k(t,x)-u_n(t,x)\big]_+\dx
\le \lim_{n\to\infty}\int_{\RR^d}\big[v_k(0,x)-u_n(0,x)\big]_+\dx
=\int_{\RR^d}\big[v_k(0,x)-u_0(x)\big]_+\dx=0
\]
since $v_k(0,x)\le u_0$ by construction. Therefore also $v_k(t,x)\le u(t,x)$ for $t>0$, so that in the limit $k\to \infty$ we obtain $v(t,x)\le u(t,x)$\,. The inequality $u\le v$ can be obtained simply by switching the roles of $u_n$ and $v_k$\,. The validity of estimates of Theorem $\ref{prop.HP.s}$ is guaranteed by the above limiting process. The comparison holds by taking the limits in inequality \eqref{final.001}, as it has been done for $\LL^1$-solutions in \cite{DPQRV2}.\qed


\section{Good fast diffusion range}\label{sect.3}

The first result of the section will be the existence of local lower bounds. In the proof we will use Lemma \ref{Lem.Opt}, which is a simple optimization lemma that we state in Appendix \ref{app.opt}\,. We recall that $m_c:=d/(d-2s)$ and $\vartheta:=1/[2s-d(1-m)]$ which is positive for $m>m_c$\,.

\begin{figure}[ht]
\centering
\ifpdf
    \includegraphics[height=8.2cm, width=15cm]{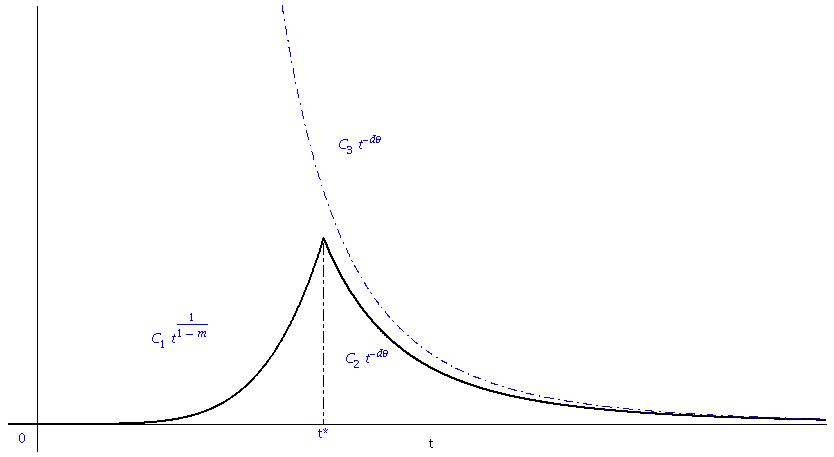}
\else
    \includegraphics[height=8.2cm, width=15cm]{timesGFDE.eps}
\fi
 \caption{\noindent\textit{Black: Lower bounds in the two time ranges. Blue: Upper bounds (smoothing effects), which has the same behaviour when $t\ge t_*$\,.}}
  \label{fig.1}
\end{figure}

\begin{thm}[Local lower bounds]\label{thm.lower}
Let $R_0>0$, $m_c<m<1$ and let $0\le u_0\in\LL^1(\RR^d, \varphi\dx)$, where $\varphi$ is as in Theorem $\ref{prop.HP.s}$ with decay at infinity $|x|^{-\alpha}$, $d-[2s/(1-m)]<\alpha<d+(2s/m)$. Let $u(t,\cdot)\in\LL^1(\RR^d, \varphi\dx)$ be a very weak solution to Equation \eqref{FDE.eq} corresponding to the initial datum $u_0$.  Then there exists a time
\begin{equation}\label{t*}
t_*:=C_* \,R_0^{2s-d(1-m)}\,\|u_0\|_{\LL^1(B_{R_0})}^{1-m}
\end{equation}
such that
\begin{equation}\label{low.1.thm}
\inf_{x\in B_{R_0/2}}u(t,x)\ge
K_1\,R_0^{-\frac{2s}{1-m}}\,t^{\frac{1}{1-m}}\quad \mbox{ if } \ 0\le t\le t_*\,,
\end{equation}
and
\begin{equation}
\inf_{x\in B_{R_0/2}}u(t,x)\ge K_2\dfrac{\|u_0\|_{\LL^1(B_{R_0})}^{2s\vartheta}}{t^{d\vartheta}}\quad \mbox{ if } \ t\ge t_*\,.
\end{equation}
The positive constants $C_*,K_1,K_2$  depend only on $m,s$ and $d\ge 1$. See Fig. 1.
\end{thm}

\medskip

\noindent\textbf{Remarks. }(i) The lower estimate for small times is an absolute bound in the sense that it does not depend on the initial data (though $t_*$ does depend).\\
(ii) We obtain the following expressions for $K_1$ and $K_2$ and $C_*$:
\begin{equation}\label{K.i}\begin{split}
&K_1:= \dfrac{K_2}{\left[2^{\frac{2}{\vartheta}+1}s\vartheta\left(\omega_d\,I_\infty\right)^{\frac{1}{d\vartheta}}
    \right]^{d\vartheta+\frac{1}{1-m}}}\,,\qquad\mbox{and}\qquad \\
&K_2:=\left[\left(\frac{2s}{d(1-m)}\right)^{\frac{1}{\vartheta}}-1\right]^{\frac{1}{1-m}}
\left[\frac{d(1-m)}{2s} \frac{\left(2s\vartheta\,\right)^{d(1-m)\vartheta}-1}{\left(2s\vartheta\,\right)^{d(1-m)\vartheta}}\right]^{\frac{2s\vartheta}{1-m}}
\frac{\alpha-d}{2(\alpha-d)+1}\,\frac{1}{\omega_d  4^d C_1^{d\vartheta}}\\
&C_* = 2s\vartheta\left(\omega_d\, 2^d\, I_\infty\right)^{\frac{1}{d\vartheta}}
\end{split}
\end{equation}
where $C_1>0$ is the constant in the $\LL^1$-weighted estimates of Proposition \ref{prop.HP.s} that depends on $\alpha,m,d$, with $d<\alpha< d+\frac{2s}{m}$, and $I_\infty>0$ is the constant in the smoothing effects \eqref{Smoothing.FFDE}, cf. Theorem 2.2 of \cite{DPQRV2}.\\
(iii) We can always choose $\alpha=d/m<d+2s/m$\,, since $2s>d(1-m)$.

\medskip

\noindent {\sl Proof.~}The proof is divided in several steps.

\noindent$\bullet~$\textsc{Step 1. }\textit{Reduction. }By the comparison principle that it is sufficient to prove lower bounds for solutions $u$ to the following reduced problem:
\begin{equation}\label{FFDE.Prob.Red}
\left\{
\begin{array}{lll}
\partial_t u
+(-\Delta)^s (u^m)=0\,,\; &\mbox{in }(0,\infty)\times\RR^d\,,\\
u(0,\cdot)=u_0\chi_{B_{R_0}}=\underline{u_0}\,,\; &\mbox{in }\RR^d\,,
\end{array}
\right.
\end{equation}
where $m_c<m<1$\,, $0<s<1$\,, and $R_0>0$\,. We only assume that $0\le u_0\in\LL^1(B_{R_0})$\,, which implies that $\underline{u_0}\in \LL^1(\RR^d)$ since $\supp(\underline{u_0})\subseteq B_{R_0}$ and also that $\|\underline{u_0}\|_{\LL^1(\RR^d)}=\|u_0\|_{\LL^1(B_{R_0})}$\,. It is not restrictive to assume that the ball $B_{R_0}$ is centered at the origin.

\noindent$\bullet~$\textsc{Step 2. }\textit{Smoothing effects. }In \cite{DPQRV2} there are the global $\LL^1-\LL^\infty$ smoothing effects which provide global upper bounds for solutions to the Cauchy problem \ref{FDE.eq}\,. We apply such smoothing effects to solutions to our reduced Problem \ref{FFDE.Prob.Red} to get
\begin{equation}\label{Smoothing.FFDE}
\|u(t)\|_{\LL^\infty(\RR^d)}\le \frac{I_\infty}{t^{d\vartheta}}\|\underline{u_0}\|_{\LL^1(\RR^d)}^{2s\vartheta}
=\frac{I_\infty}{t^{d\vartheta}}\|u_0\|_{\LL^1(B_{R_0})}^{2s\vartheta}
\end{equation}
where $\vartheta=1/[2s-d(1-m)]$ and the constant $I_\infty$ only depends on $d,s,m$\,.

\noindent$\bullet~$\textsc{Step 3. }\textit{Aleksandrov principle. }We recall Theorem 11.2 of \cite{Vaz2012}, we have that
\[
u(t,0)\ge u(t,x)\,,\qquad\mbox{for all}\; t>0\;\mbox{and}\; |x|\ge 2R_0\,.
\]
Therefore one has that
\begin{equation}\label{aleks.1}
\|u(t)\|_{\LL^\infty(\RR^d\setminus B_{2R_0})}=\sup_{x\in \RR^d\setminus B_{2R_0}}u(t,x)\le u(t,0)\,.
\end{equation}

\noindent$\bullet~$\textsc{Step 4. }\textit{Lower estimates for the $\LL^\infty$-norm on an annulus. }We combine the $\LL^1$-weighted estimates of Theorem \ref{prop.HP.s} with the smoothing effects of Step 2: estimates \eqref{HP.s} read in this context
\begin{equation}\label{HP.s.step41}
\left(\int_{B_{R_0}}u_0\dx\right)^{1-m}\le \left(\int_{\RR^d}u_0\varphi_R(x)\dx\right)^{1-m}\le
\left(\int_{\RR^d}u(t,x)\varphi_R(x)\dx\right)^{1-m}
+ C_1 R^{d(1-m)-2s}\,t
\end{equation}
we have chosen $R\ge 2 R_0>0$ and $\varphi_R(x)=\varphi(x/R)$ with $\varphi$ as in Lemma \ref{Lem.phi}  (with the explicit form given in formula \eqref{phi})\,,  so that $\varphi_R(x)=1$ on $B_R$ and $0\le \varphi_R(x)\le |x|^{-\alpha}$ for $|x|>>R$ with $d-2s/(1-m)<\alpha< d+2s/m$\,, and we recall that $C_1>0$ depends only on $\alpha,m,d$\,.
\begin{equation}\label{HP.s.2}\begin{split}
&\|u_0\|_{\LL^1(B_{R_0})}^{1-m}
- C_1 R^{d(1-m)-2s}\,t
 \le \left(\int_{\RR^d}u(t,x)\varphi_R(x)\dx\right)^{1-m}\\
&\le \left(\int_{\RR^d\setminus B_{2 R_0}}u(t,x)\varphi_R(x)\dx\right)^{1-m}
    + \left(\int_{B_{2 R_0}}u(t,x)\varphi_R(x)\dx\right)^{1-m}=(I)+(II)\,.\\
\end{split}
\end{equation}
We first estimate (I), to this end we observe that if we choose $d<\alpha< d+2s/m$ we have that
\begin{equation}\label{ineq.phi.R}
\begin{split}
\int_{\RR^d\setminus B_{2 R_0}}\varphi_R(x)\dx
    &=\int_{\RR^d\setminus B_{R}}\varphi_R(x)\dx +\int_{B_R\setminus B_{2 R_0}}\varphi_R(x)\dx
     =\int_{\RR^d\setminus B_{R}}\varphi_R(x)\dx +\int_{B_R\setminus B_{2 R_0}}1\dx\\
    &\le\int_{\RR^d\setminus B_{R}}\dfrac{1}{\big[1+(|x/R|^2-1)^4\big]^{\alpha/8}}\dx+\omega_d R^d\\
    &=\omega_d R^d \int_{1}^{+\infty}\dfrac{r^{d-1}}{\big[1+(r^2-1)^4\big]^{\alpha/8}}\rd r +\omega_d R^d\\
    &=\omega_d R^d \left[\int_{1}^{4}\dfrac{r^{d-1}}{\big[1+(r^2-1)^4\big]^{\alpha/8}}\rd r
     +\int_{4}^{+\infty}\dfrac{r^{d-1}}{\big[1+(r^2-1)^4\big]^{\alpha/8}}\rd r +1\right]\\
    &\le_{(a)}\omega_d R^d \left[1+ 4^d+4^{\alpha/8}\int_{4}^{+\infty}r^{d-1-\alpha}\rd r\right]
    =\omega_d R^d \left[1+4^d+\frac{4^{\alpha/8}}{\alpha-d}\frac{1}{4^{\alpha-d}}\right]\\
    &\le \omega_d  4^d \frac{2(\alpha-d)+1}{\alpha-d}R^d\\
\end{split}
\end{equation}
where we have used that $R\ge 2 R_0$ and in $(a)$ we have used the fact that $\varphi_R\le 1$ and that  $1+\big(r^2-1\big)^4\ge r^8/4$\,, if $r\ge 4$\,. Therefore we have
\[
\begin{split}
(I)^{\frac{1}{1-m}}
&=\int_{\RR^d\setminus B_{2 R_0}}u(t,x)\varphi_R(x)\dx
\le \|u(t)\|_{\LL^\infty(\RR^d\setminus B_{2 R_0})}\int_{\RR^d\setminus B_{2 R_0}}\varphi_R(x)\dx\\
&\le  \omega_d  4^d \frac{2(\alpha-d)+1}{\alpha-d}R^d \|u(t)\|_{\LL^\infty(\RR^d\setminus B_{2 R_0})}
\le  \omega_d  4^d \frac{2(\alpha-d)+1}{\alpha-d}\,R^d\, u(t,0)
\end{split}
\]
where in the last step we have used inequality \eqref{aleks.1} of Step 3, derived from Aleksandrov principle.

\noindent We now estimate $(II)$ as follows:
\[\begin{split}
(II)^{\frac{1}{1-m}}
    &=\int_{B_{2 R_0}}u(t,x)\varphi_R(x)\dx\le \|u(t)\|_{\LL^\infty(\RR^d)}\int_{B_{2 R_0}}\varphi_R(x)\dx\\
    &\le_{(a)} \omega_d\, 2^d \,R_0^d\, \|u(t)\|_{\LL^\infty(\RR^d)}
     \le_{(b)} \omega_d\, 2^d\, R_0^d \,\frac{I_\infty}{t^{d\vartheta}}\,\|u_0\|_{\LL^1(B_{R_0})}^{2s\vartheta}
\end{split}\]
where in $(a)$ we have used that $\varphi_R(x)=1$ on $B_R$, $2 R_0<R$ and $|B_R|=\omega_d R^d$\,. In $(b)$ we have used the smoothing effect \eqref{Smoothing.FFDE}\,. Plugging the above estimates into \eqref{HP.s.2} gives
\begin{equation}\label{HP.s.33}\begin{split}
\|u_0\|_{\LL^1(B_{R_0})}^{1-m} - C_1 R^{d(1-m)-2s}\,t
&\le   \,\frac{\left[\omega_d\, 2^d\, R_0^d I_\infty\,\|u_0\|_{\LL^1(B_{R_0})}^{2s\vartheta}\right]^{1-m}}{t^{d(1-m)\vartheta}}\,\\
    &+\left[\omega_d  4^d \frac{2(\alpha-d)+1}{\alpha-d}\right]^{1-m}\,R^{d(1-m)}\, u^{1-m}(t,0) \,,\\
\end{split}
\end{equation}
or equivalently
\begin{equation}\label{HP.s.4}\begin{split}
\left[ \|u_0\|_{\LL^1(B_{R_0})}^{1-m}
  -\,\frac{\left[\omega_d\, 2^d\, R_0^d I_\infty\,\|u_0\|_{\LL^1(B_{R_0})}^{2s\vartheta}\right]^{1-m}}{t^{d(1-m)\vartheta}}\right]&\frac{1}{R^{d(1-m)}}\,
    - \frac{C_1\, t}{R^{2s}}\\
    &\le \left[\omega_d  4^d \frac{2(\alpha-d)+1}{\alpha-d}\right]^{1-m}\, u^{1-m}(t,0)\,.\\
\end{split}
\end{equation}

\noindent$\bullet~$\textsc{Step 5. }\textit{Optimization. }The previous estimate \eqref{HP.s.4} is useful only if we can make sure that the left-hand side has a positive lower bound.  Let us write inequality \eqref{HP.s.4} as
\begin{equation}\label{HP.s.5}\begin{split}
F(t,R):=&\frac{A(t)}{R^{d(1-m)}}\,
    - \frac{B\, t}{R^{2s}}
    \le \left[\omega_d  4^d \frac{2(\alpha-d)+1}{\alpha-d}\right]^{1-m}\, u^{1-m}(t,0)\,,\\
\end{split}
\end{equation}
with
\begin{equation}\label{F.A.B.C}\begin{split}
A(t)=\left[M-\,\frac{C}{t^{d(1-m)\vartheta}}\right]\,,\; M:=\|u_0\|_{\LL^1(B_{R_0})}^{1-m}\,,\; C:=\left[\omega_d\, 2^d\, R_0^d I_\infty\,\|u_0\|_{\LL^1(B_{R_0})}^{2s\vartheta}\right]^{1-m}\,,\; B=C_1
\end{split}
\end{equation}
where $C_1>0$ is the constant of $\LL^1$-weighted estimates of Theorem \ref{prop.HP.s}, and $I_\infty>0$ is the constant of the smoothing effects \eqref{Smoothing.FFDE} of Step 2. We now optimize the function $F$ as in Lemma \ref{Lem.Opt} so that
there exists
\begin{equation}\label{t.min.step5}
t_*:=2s\vartheta\left(\frac{C}{M}\right)^{\frac{1}{d(1-m)\vartheta}}
=2s\vartheta\left(\omega_d\, 2^d\, I_\infty\right)^{\frac{1}{d\vartheta}}\,R_0^{\frac{1}{\vartheta}}\,\|u_0\|_{\LL^1(B_{R_0})}^{1-m}
\end{equation}
and
\begin{equation}\label{R.max.step5}
\begin{split}
\overline{R}(t)
&=\left(\frac{2sBt}{d(1-m)A(t)}\right)^{\vartheta}
    \ge\overline{R}(t_*)=\left[\frac{2s}{d(1-m)}\frac{(2s\vartheta)^{2s\vartheta}}{(2s\vartheta)^{d(1-m)}-1}\right]^{\vartheta}
    \frac{B^\vartheta C^{\frac{1}{d(1-m)}}}{M^{\frac{2s\vartheta}{d(1-m)}}}\\
&=\left[\frac{2s}{d(1-m)}\frac{(2s\vartheta)^{2s\vartheta}}{(2s\vartheta)^{d(1-m)}-1}\right]^{\vartheta}
 \omega_d^{\frac{1}{d}}\, 2\, R_0 I_\infty^{\frac{1}{d}}C_1^\vartheta \ge_{(a)}2 R_0\,,
\end{split}
\end{equation}
where in $(a)$ we have used that the constants $I_\infty>0$ and $C_1>0$ are constants in the upper bounds \eqref{Smoothing.FFDE} and \eqref{HP.s.step41} respectively, so that we can chose them to be arbitrarily large to fulfill the condition $\overline{R}(t_*)\ge 2 R$. Therefore for all $t\ge t_*$ we have that
\[\begin{split}
\left[\omega_d  4^d \frac{2(\alpha-d)+1}{\alpha-d}\right]^{1-m}\, u^{1-m}(t,0)
&\ge F(\overline{R}(t),t)
=\left[\left(\frac{2s}{d(1-m)}\right)^{\frac{1}{\vartheta}}-1\right]\left[\frac{d(1-m)}{2s}\right]^{2s\vartheta}
\frac{A(t)^{2s\vartheta}}{(Bt)^{d(1-m)\vartheta}}\\
&\ge \left[\left(\frac{2s}{d(1-m)}\right)^{\frac{1}{\vartheta}}-1\right]\left[\frac{d(1-m)}{2s}\right]^{2s\vartheta}
\frac{A(t_*)^{2s\vartheta}}{C_1^{d(1-m)\vartheta}}\frac{1}{t^{d(1-m)\vartheta}}
\end{split}\]
since $A(t)\ge A(t_*)$ for all $t\ge t_*$, and it is easy to check that
\[
A(t_*)=\left[1-\,\frac{1}{\left(2s\vartheta\,\right)^{d(1-m)\vartheta}}\right]
\|u_0\|_{\LL^1(B_{R_0})}^{1-m}
=\frac{\left(2s\vartheta\,\right)^{d(1-m)\vartheta}-1}{\left(2s\vartheta\,\right)^{d(1-m)\vartheta}}
\|u_0\|_{\LL^1(B_{R_0})}^{1-m}>0\,,\;
\]
since we recall that $2s\vartheta>1$. Summing up we have obtained
\begin{equation}\label{step5.1}
u(t,0)\ge K_2 \frac{\|u_0\|_{\LL^1(B_{R_0})}^{2s\vartheta}}{t^{d\theta}}\,,
\end{equation}
for all $t\ge t_*>0$\,, where $K_2$ only depends from $\alpha,m,s,d$ and takes the form
\begin{equation}\label{step5.1.K}
K_2:=\left[\left(\frac{2s}{d(1-m)}\right)^{\frac{1}{\vartheta}}-1\right]^{\frac{1}{1-m}}
\left[\frac{d(1-m)}{2s} \frac{\left(2s\vartheta\,\right)^{d(1-m)\vartheta}-1}{\left(2s\vartheta\,\right)^{d(1-m)\vartheta}}\right]^{\frac{2s\vartheta}{1-m}}
\frac{\alpha-d}{2(\alpha-d)+1}\,\frac{1}{\omega_d  4^d C_1^{d\vartheta}}
\end{equation}
Note that in the limit $m\to 1$ the constant $K_2\to 0$.   By a standard argument it is easy to pass from the center to the infimum on $B_{R_0/2}(0)$ in the above estimates.

\noindent$\bullet~$\textsc{Step 6. }\textit{Positivity backward in time. }Using Benilan-Crandall estimates which depend only by the homogeneity of the equations, cf. \cite{BCr}
\begin{equation}\label{BC.est}
u_t\le \frac{u}{(1-m)t}
\end{equation}
we can prove positivity in the time interval $[0,t_*]$. These estimates in the fractional case has been proven in \cite{DPQRV2}, and imply
that the function: $u(t,x)t^{-1/(1-m)}$ is non-increasing in time, thus for any $t\in (0,t_*)$ and $x\in B_{R_0/2}(0)$\,, inequality \eqref{step5.1} gives
\[
u(t,x)\geq \frac{u(t_*,x)}{t_*^{\frac{1}{1-m}}}t^{\frac{1}{1-m}}\ge K_2 \frac{\|u_0\|_{\LL^1(B_{R_0})}^{2s\vartheta}}{t_*^{d\theta+\frac{1}{1-m}}}t^{\frac{1}{1-m}}
=\frac{K_2}
    {\left[2^{\frac{2}{\vartheta}+1}s\vartheta\left(\omega_d\,I_\infty\right)^{\frac{1}{d\vartheta}}
    \right]^{d\theta+\frac{1}{1-m}}}\left[\frac{t}{R_0^{2s}}\right]^{\frac{1}{1-m}}
\]
where $K_2>0$ is given in \eqref{step5.1.K}\,, and $t_*$ is given by \eqref{t.min.step5}, and it is easy to check that
\begin{equation*}\label{step7.1}\begin{split}
\frac{\|u_0\|_{\LL^1(B_{R_0})}^{2s\vartheta}}{t_*^{d\theta+\frac{1}{1-m}}}
&=\frac{\|u_0\|_{\LL^1(B_{R_0})}^{2s\vartheta}}
    {\left[2s\vartheta\left(\omega_d\, 2^d\,I_\infty\right)^{\frac{1}{d\vartheta}}\,(R_0)^{\frac{1}{\vartheta}}
    \|u_0\|_{\LL^1(B_{R_0})}^{1-m}\right]^{d\theta+\frac{1}{1-m}}}
=\frac{1}
    {\left[2^{\frac{2}{\vartheta}+1}s\vartheta\left(\omega_d\,I_\infty\right)^{\frac{1}{d\vartheta}}
    \right]^{d\theta+\frac{1}{1-m}}R_0^{\frac{2s}{1-m}}}\,.
\end{split}
\end{equation*}
The proof is concluded.\qed

\medskip

{\bf  Remark.} This lower estimate holds in the limit $m\to 1$ and gives lower estimates for the linear fractional heat equation of the following form.

\begin{prop}\label{thm.lower.m=1}
Let $u\ge 0$ be a weak solution to the Cauchy Problem \eqref{FDE.eq}, corresponding to $u_0\in \LL^1(\RR^d)$ and $m=1$. Then  $\vartheta=1/2s>0$ and the estimate says that for given $R_0$ and
$t_*:=C_1 \,R_0^{2s}$, then
\begin{equation}
\inf_{x\in B_{R_0/2}}u(t,x)\ge
K_2\dfrac{\|u_0\|_{\LL^1(B_{R_0})}}{t^{d/2s}}\quad \mbox{ if } \ t\ge t_*\,.
\end{equation}
The positive constant $K_2$  depends only on $C_1$, $s$ and $d$.
\end{prop}

The proof is easily obtained from the integral representation of the solution.

\subsection{Minimal space-like tail behaviour}\label{ssec.gfd.tail}

As a corollary of the previous lower bound, we obtain a quantitative bound from below for the space-like behaviour of any nonnegative solution. We consider a solution that has a certain initial mass $M$ in the ball of radius 1 and apply the result of Theorem \ref{thm.lower}  after displacing the origin of space coordinates to a point $x_0$ with $|x_0|>>1$. We then consider the formula \eqref{t*} for the critical time with center $x_0$ and radius $R_0=|x_0|+2$, so that the ball $B_{R_0}(x_0)$ contains the mass $M$ mentioned above. As $R_0\to \infty$ also $t_*\to\infty$. We can therefore use the lower bound \eqref{low.1.thm} to get an estimate of the form
\begin{equation}\label{low.bdd.20}
u(t,x_0)\ge G(u_0,t)\,|x_0|^{-2s/(1-m)}\,,
\end{equation}
where $G(u_0,t)$ is given in \eqref{low.1.thm}.
According to the results of \cite{Vaz2012} the Barenblatt solutions have this precise spatial behaviour in the range $m_c<m<m_1$, with $m_1=d/(d+2s)$, therefore the asymptotic estimate is sharp in this range.

\noindent \subsection{Global spatial lower bounds in the case $m_1<m<1$}\label{ssec.global.GFDE}

We would like to prove that the solution can always be bounded from below by a Barenblatt solution, so the lower bound will be sharp.
\begin{figure}[ht]
\centering
\ifpdf
    \includegraphics[height=7cm, width=9.8cm]{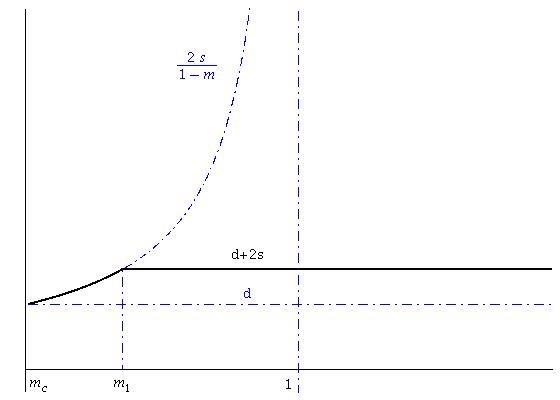}
\else
    \includegraphics[height=7cm, width=9.8cm]{tailsGFDE.eps}
\fi
 \caption{\noindent\textit{Lower bounds for the spatial decay rates of solutions. Recall that $m_c=(d-2s)/d$ and $m_1=d/(d+2s)$}}
 \label{fig.2}
\end{figure}
In the range $m_c<m<m_1$ the lower bound of Theorem $\ref{thm.lower}$, gives sharp lower bounds with the same tails as the Barenblatt solutions, as explained in Section \ref{ssec.gfd.tail}. In the range $m_1<m\le 1$ the lower bound given by \eqref{low.bdd.20} is not sharp and the following Theorem \ref{thm.lower.m1} (respectively Proposition $\ref{thm.lower.m=1}$ when $m=1$) proves that any solution with data in $\LL^1(\RR^d)$ can always be bounded from below by a Barenblatt solution (respectively by the fundamental solution when $m=1$). See Figure \ref{fig.2}.

\begin{thm}[Global Lower Bounds when $m_1<m<1$]\label{thm.lower.m1} Under the conditions of Theorem $\ref{thm.lower}$ we have in the range  $m_1<m<1$
\begin{equation}
u(t,x)\ge \frac{C(t)}{|x|^{d+2s}}\qquad \mbox{ when $|x|>>1$.}
\end{equation}
valid for all $0<t<T$ with some bounded function $C>0$ that depends on $t,T$ and on the data.
\end{thm}

\noindent {\sl Proof.~}The proof consists of several steps.

\noindent$\bullet~$\textsc{Step 1. } We begin under the extra assumption that $u_0(x)\ge 2c>0$ in a ball, that can be taken to be $B_1(0)$ by scaling.  Therefore, there exists a $t_1$ such that
$u(t,x)\ge c$ for all $0\le t\le t_1$ and all $|x|\le 1$. We also assume that $u_0$ is continuous and goes to zero uniformly as $|x|\to\infty$.

Consider the function $u_{0,\varepsilon}(x)=u_0(x)+\varepsilon$\,, and let $u_\varepsilon(t,x)\ge u(t,x)$ be the corresponding solutions. By the usual theory, \cite{DPQRV2}, we know that $u_\varepsilon\ge \varepsilon$, $u_\varepsilon-\varepsilon\in \LL^1(\RR^d)$, since $u_0\in \LL^1(\RR^d)$.  Moreover, it is proved in the theory that $u_\varepsilon\to \varepsilon$ as $|x|\to\infty$ for every $t>0$.

\noindent$\bullet~$\textsc{Step 2. }Let $\underline{u}=u^*(t+\tau,x)$, where $u^*$ is the Barenblatt solution with mass $M>0$. We refer to \cite{Vaz2012} for a complete discussion about Barenblatt solutions.  By choosing the mass $M:=\int_{\RR^d}\underline{u}\dx$ very small, we can find $\tau=\tau(\varepsilon)>0$ so that $\underline{u}(t,x)\le c/4$ for $|x|\ge \delta$ and $0<t<t_1$, and $\underline{u}\le \varepsilon/2$ for $|x|\ge 1$.

\noindent$\bullet~$\textsc{Step 3. }We compare both continuous solutions in the exterior domain $\Omega=\{x: |x|\ge 1\}$. At the first time where $\underline{u}$ touches $u_\varepsilon$ from below at a point $|x|> 1$, we have $\partial_t(u_\varepsilon-\underline{u})\le 0$. Let now $w=u_\varepsilon^m-\underline{u}^m$ to get
\[\begin{split}
(-\Delta)^s w(x) &=k_{s,d}\int_{\RR^d}\frac{w(x)-w(y)}{|x-y|^{d+2s}}\dy \\
&= k_{s,d}\int_{\{|x|\le 1\}}\frac{w(x)-w(y)}{|x-y|^{d+2s}}\dy
+k_{s,d}\int_{\{|x|\ge 1\}}\frac{w(x)-w(y)}{|x-y|^{d+2s}}\dy=I_1+I_2
\end{split}
\]
We want to prove now that both $I_1< 0$ and $I_2< 0$, which leads to a contradiction. In this way we conclude that $\underline{u}<u_\varepsilon$ for all $0<t<t_1$ and all $|x|\ge 1$.

The fact that $I_2\le 0$ comes from the fact that $w(x)=0$ by our choice of $x$ and $w(y)\ge 0$ since $(x,t)$ is the first contact point. Due to the fact that $w(y)>0$ near $|x|=1$ and $w$ is continous we get $I_2<0$.

As for $I_1$, the denominator is like a constant in that domain and we have to estimate $w(y)$. We know that for $\delta<|x|<1$ we have $u_\varepsilon \ge c$ and $\underline{u}\le c/4$, hence $w(y)\ge C\,c^{m}>0$ and this contributes to the integral something that is like $-C\,c^{m}$, which is not small. In the small ball $|x|\le \delta $ we use the worst case estimate
$-w(y)\le \underline{u}$ and $\underline{u}(t_1,y)$ has mass at most $M$ which is small, this contributes
at most a bad term of order
$$
C\int_{|x|\le \delta} \underline{u}^m\,dx\le C M^m \,\delta^{d(1-m)},
$$
which is small if $\delta $ and $M$ are small (here we use $m<1$).  Therefore $I_1<0$.

Moreover, one has to ensure that $\underline{u}(0,x)<u_\varepsilon(0,x)$ for $|x|\ge 1$. Since $u_\varepsilon\ge \varepsilon$ and
\[
\underline{u}(0,x)=\tau^{-\alpha} F(|x|\,\tau^{-\beta})\le \tau^{-\alpha} F(\tau^{-\beta})=c \tau^{-\alpha+\beta(d+2s)}=c\tau^{2s\beta}\le \varepsilon
\]
at least for sufficiently small $\tau$, recall Step 2. By the parabolic comparison theorem we conclude that  $\underline{u}<u_\varepsilon$ for all $0<t<t_1$ and all $|x|\ge 1$.

\medskip

\noindent$\bullet~$\textsc{Step 4. }We finally let $\varepsilon\to 0$ and also $\tau$ may go to zero, and we obtain that $u^*(x,t):=\lim_{\varepsilon\to 0}u_\varepsilon= u$, therefore we can conclude that $u(x,t)\ge c/|x|^{d+2s}$ when $|x|>>1$ and $t=t_1$.

\medskip

\noindent$\bullet~$\textsc{Step 5. }Once we have obtained the spatial lower bound at times $t\le t_1$, then we can compare with a Barenblatt solution and continue the lower bound for all times, to finally get that the spatial tail of the solution $u$ can be bounded from below by $u\ge c/|x|^{d+2s}$ when $|x|>>1$.\qed

\section{Very fast diffusion range}\label{sect.4}

In the very fast diffusion range $0<m<m_c$, the weighted $\LL^1$ estimates of Theorem \ref{prop.HP.s} continue to hold, but this does not allow to obtain quantitative lower bounds since technique used in the good fast diffusion range does not work anymore. One problem is that the smoothing effect does not hold for general $\LL^1$ initial data, therefore the optimization of Lemma \ref{Lem.Opt} is no more valid, since $2s<d(1-m)$ in this range. Hence the need for new weighted $\LL^1$ estimates, in the form given in Step 3 of the proof of Theorem \ref{thm.lower.subcrit.} below. Another problem typical of this range of exponent is the presence of the extinction time, which enters directly in the estimates of Theorem \ref{thm.lower.subcrit.}. We present here a technique that is based on the careful use of weight factors.

\begin{thm}[Local lower bounds I]\label{thm.lower.subcrit.}
Let $u$ be a weak solution to the equation \eqref{FDE.eq}, corresponding to $u_0\in \LL^1(\RR^d)\cap\LL^{p_c}(\RR^d)$ with $0<m<m_c=d/(d-2s)$, $0<s<1$ and let $p_c=d(1-m)/(2s)$. Let also $T=T(u_0)$ be the finite extinction time for $u$. Then for every $R_0>0$, there exists a time
\begin{equation}\label{t*.subcrit}
t_*:= C_*\,R_0^{2s-d(1-m)}\|u_0\|_{\LL^1(B_{R_0})}^{1-m}\le T(u_0)\,,
\end{equation}
such that
\begin{equation}\label{thm.pos.T.VFDE}
\inf_{x\in B_{R_0/2}}u(t,x)\ge K \frac{\|u_0\|_{\LL^1(B_{R_0})}^{\frac{1}{m}} }{R_0^{\frac{d-2s}{m}} } \frac{t^{\frac{1}{1-m}}}{T^{\frac{1}{m(1-m)}}}
\quad \mbox{ if } \ 0\le t\le t_*\,,
\end{equation}
where $C_*$ and $K$ are explicit positive universal constants, that depend only on $m,s,d$.
\end{thm}
The expression of the constants is
\begin{equation}\label{const.subcrit}\begin{split}
 C_*&:=\frac{k_{s,d}\,\omega_d^m}{4^{d+1-2s}}\,,\qquad
K :=\left(\frac{k_{s,d}}{4^{3d+1-2s}\,d}\right)^{\frac{1}{m}}\,,\\
\end{split}\end{equation}\normalcolor
where $k_{s,d}$ is the constant of the representation formula $  \varphi(x)=k_{s,d}\int_{\RR^d}\frac{\rho(y)}{|x-y|^{d-2s}}\dy$ and $\omega_d$ is the volume of the unit ball.

\noindent {\sl Proof. of Theorem \ref{thm.lower.subcrit.}~} It is divided into several steps as follows.

\noindent$\bullet~$\textsc{Step 1. }\textit{Reduction. }By the comparison principle that it is sufficient to prove lower bounds for solutions $u$ to the following reduced problem:
\begin{equation}\label{FDE.eq.Red.VFDE}
\left\{
\begin{array}{lll}
\partial_t u
+(-\Delta)^s (u^m)=0\,,\; &\mbox{in }(0,\infty)\times\RR^d\,,\\
u(0,\cdot)=u_0\chi_{B_{R_0}}=\underline{u_0}\,,\; &\mbox{in }\RR^d\,,
\end{array}
\right.
\end{equation}
where $m>1$\,, $0<s<1$\,, and $R_0>0$\,. We only assume that $0\le u_0\in\LL^1(B_{R_0})$\,, which implies that $\underline{u_0}\in \LL^1(\RR^d)$ since $\supp(\underline{u_0})\subseteq B_{R_0}$ and also that $\|\underline{u_0}\|_{\LL^1(\RR^d)}=\|u_0\|_{\LL^1(B_{R_0})}$\,. It is not restrictive to assume that the ball $B_{R_0}$ is centered at the origin. We call $M_0=\|u_0\|_{\LL^1(B_{R_0})}$.

\noindent$\bullet~$\textsc{Step 2. }\textit{Aleksandrov principle. } We recall Theorem 11.2 of \cite{Vaz2012}. In view of the fact that the initial function is supported in the ball $B_{R_0}(0)$,  we have that
\[
u(t,0)\ge u(t,x)\,,\qquad\mbox{for all}\; t>0\;\mbox{and}\; |x|\ge 2R_0\,.
\]
Therefore, one has
\begin{equation}\label{aleks.1.PME}
\sup_{x\in \RR^d\setminus B_{2R_0}}u(t,x)\le u(t,0)\,.
\end{equation}

\noindent$\bullet~$\textsc{Step 3. }\textit{$\LL^1$ Weighted estimates. }Choose a test function $\varphi\ge 0$ such that $-(-\Delta)^s\varphi=\rho$ with $\rho=0$ on $B_{2R_0}$ and on $B_{R_1}^c$, and $0< \rho \le  1$ on the annulus $A:= B_{R_1}\setminus B_{2R_0}$, with $0<2R_0\le R_1$, and $R_0$ as in Step 1, such that $\supp(u_0)\subseteq B_{R_0}$.
Using the explicit representation of $\varphi$ in terms $\rho$ and the integral kernel $K(x,y)=k_{s,d}|x-y|^{n-2s}$ we get the estimates
\[
\varphi(x)=k_{s,d}\int_{\RR^d}\frac{\rho(y)}{|x-y|^{d-2s}}\dy\ge \frac{k_{s,d} \|\rho\|_1}{(R_1+R_0)^{d-2s} }\ge k_0>0\,,\qquad\mbox{for all $x\in B_{R_0}(0)$\,,}
\]
since $|x-y|\le R_0+R_1$. We can always choose $\rho\ge 1/2$ on the smaller annulus $A_0= B_{2R_0+3(R_1-2R_0)/4}\setminus B_{2R_0+(R_1-2R_0)/4}\subseteq A$\,, so that
\[
\begin{split}
\|\rho\|_1
    &= \int_{A_1}\rho(x)\dx
     \ge\int_{A_0}\rho(x)\dx
     \ge \frac{|A_0|}{2}
     =\frac{1}{2}|B_{2R_0+3(R_1-2R_0)/4}\setminus B_{2R_0+(R_1-2R_0)/4}|\\
     &=\frac{\omega_d}{2}\left[\left(2R_0+\frac{1}{4}(R_1-2R_0) + (R_1-2R_0)\right)^d-\left(2R_0+\frac{1}{4}(R_1-2R_0)\right)^d\right]
     \ge \frac{\omega_d}{2}(R_1-2R_0)^d
\end{split}
\]
since $(a+b)^d- a^d\ge b^d$ for any $a,b\ge 0$.  Then $k_0>0$ takes the form
\begin{equation}\label{k0.VFFDE}
k_0:=\frac{k_{s,d}\,\omega_d}{2}\frac{(R_1-2R_0)^d}{(R_1+R_0)^{d-2s}}
\end{equation}\normalcolor
Now we observe that letting $T=T(u_0)>0$ be the finite extinction time for the reduced problem \eqref{FDE.eq.Red.VFDE}, we obtain
\begin{equation}\label{step.3.1}\begin{split}
\int_{\RR^d}u_0(x)\varphi\dx &=\int_{\RR^d}u(0,x)\varphi\dx - \int_{\RR^d}u(T,x)\varphi\dx
 =-\int_0^T\int_{\RR^d}\partial_\tau u(\tau,x)\varphi(x)\dx \rd\tau\\
 &= \int_0^T\int_{\RR^d}(-\Delta)^s \big(u^m(\tau,x)\big)\,\varphi(x)\dx \rd\tau
 = \int_0^T\int_{\RR^d}u^m(\tau,x)(-\Delta)^s\varphi(x)\dx \rd\tau\\
 &= \int_0^T\int_{A}u^m(\tau,x)\rho(x) \dx \rd\tau\\
 &= \int_0^{t_*}\int_{A}u^m(\tau,x)\rho(x)\dx \rd\tau + \int_{t_*}^T\int_{A}u^m(\tau,x)\rho(x)\dx \rd\tau
:=(I)+(II)\\
\end{split}
\end{equation}
where $0\le t_*\le T$ will be chosen later.

Next we estimate $(I)$. We first observe that
\begin{equation*}
\begin{split}
\int_{A}u^m(\tau,x)\rho(x)\dx
&\le  \int_{B_{R_1}}u^m(\tau,x)\dx
\le  |B_{R_1}|^{1-m}\left(\int_{B_{R_1}}u(\tau,x)\dx\right)^m\\
&\le  |B_{R_1}|^{1-m}\left(\int_{\RR^d}u(\tau,x)\dx\right)^m
\le  |B_{R_1}|^{1-m}\left(\int_{\RR^d}\underline{u_0}(x)\dx\right)^m\\
&\le  |B_{R_1}|^{1-m}M_0^m
\end{split}
\end{equation*}
since $0<m<1$\,, $0<\rho\le 1$, and in the last step we have used the fact that $\|u(t)\|_{\LL^1(\RR^d)}\le \|\underline{u_0}\|_{\LL^1(\RR^d)}$, which has been proven in \cite{DPQRV2}\,, together with the fact that $M_0=\|u_0\|_{\LL^1(B_{R_0})}=\|\underline{u_0}\|_{\LL^1(\RR^d)}$. Therefore,
\begin{equation}\label{step.3.2}
\begin{split}
(I):= \int_0^{t_*}\int_{A}u^m(\tau,x)\rho(x)\dx\rd\tau
&\le |B_{R_1}|^{1-m}\,t_*\,M_0^m\,.
\end{split}
\end{equation}
We now estimate (II) by using the Aleksandrov principle:
\begin{equation}\label{step.3.3}
\begin{split}
(II):= \int_{t_*}^T \int_{A}u^m(\tau,x)\rho(x)\dx\rd\tau
& \le   \int_{t_*}^T \int_{A}u^m(\tau,x)\dx\rd\tau
=_{(a)}  (T-t_*) \int_{A}u^m(\tau_1,x)\dx\\
&\le   (T-t_*) |A| \sup_{x\in A} u^m(\tau_1,x)
\le   (T-t_*) |A| u^m(\tau_1,0)
\end{split}
\end{equation}
where in $(a)$ we have used the mean value theorem for the function $U(\tau)=\int_{A}u^m(\tau,x)\dx$ so that there exists a $\tau_1\in[t_*,T]$ such that $\int_{t_*}^T U(\tau)\rd\tau=(T-t_*)U(\tau_1)$. In $(b)$ we have used the Aleksandrov principle, which gives $\sup\limits_{x\in A} u^m(\tau_1,x)\le u^m(\tau_1,0)$. Summing up, we have obtained, joining \eqref{step.3.1}, \eqref{step.3.2} and \eqref{step.3.3}
\begin{equation}\label{step.3.4}\begin{split}
\int_{\RR^d}u_0(x)\varphi\dx
&\le |B_{R_1}|^{1-m}\,t_*\,M_0^m+ (T-t_*) |A| u^m(\tau_1,0)
\end{split}
\end{equation}
for some $\tau_1\in[t_*,T]$. In addition, we have $\int_{\RR^d}u_0(x)\varphi\dx \ge M_0 k_0$. We finally remark that from inequality \eqref{step.3.4} we get a lower bound for the extinction time, just by letting $t_*=T$ in formula \eqref{step.3.4}:
\begin{equation}\label{low.bdd.T}
k_0 M_0\le \int_{\RR^d}u_0(x)\varphi\dx
\le |B_{R_1}|^{1-m}\,T\,M_0^m\,,\qquad\quad\mbox{that is }\qquad\quad T\ge k_0 \frac{M_0^{1-m}}{|B_{R_1}|^{1-m}}
\end{equation}

\noindent$\bullet~$\textsc{Step 4. }\textit{Choosing the critical time $t_*$. }We now choose $t_*$ to be small enough, more precisely
\begin{equation}\label{step.4.t*}
t_*:= \frac{k_0}{2}\frac{M_0^{1-m}}{|B_{R_1}|^{1-m}}\le T\,,
\end{equation}
we note that $t_*\le T$ follows by \eqref{low.bdd.T}. With this choice of $t_*$, inequality \eqref{step.3.4} becomes
\begin{equation}
\begin{split}
\frac{k_0}{2} M_0 = k_0 M_0-|B_{R_1}|^{1-m}\,t_*\,M_0^m
&\le  (T-t_*) |A| u^m(\tau_1,0)\le  T\, |A|\, u^m(\tau_1,0)
\end{split}
\end{equation}
which is the desired positivity estimate at a time $\tau_1\in [t_*,T]$\,, namely
\begin{equation}\label{step.4.1}
\begin{split}
\frac{k_0\, M_0 }{2T\, |B_{R_1}\setminus B_{2R_0}|} \le  u^m(\tau_1,0)
\end{split}
\end{equation}

\noindent$\bullet~$\textsc{Step 5. }\textit{Positivity backward in time. }Using Benilan-Crandall estimates which depend only by the homogeneity of the equations, cf. \cite{BCr}
\begin{equation}\label{BC.est.1}
u_t\le \frac{u}{(1-m)t}
\end{equation}
we can prove positivity in the time interval $[0,\tau_1]$. These estimates in the fractional case has been proven in \cite{DPQRV2}, and imply
that the function: $u(t,x)t^{-1/(1-m)}$ is non-increasing in time, thus for any $t\in [0,\tau_1]$ we have
that
\begin{equation}\label{111}\begin{split}
u(t,0)
&\ge \frac{t^{\frac{1}{1-m}}}{\tau_1^{\frac{1}{1-m}}}u(\tau_1,0)\ge
\frac{t^{\frac{1}{1-m}}}{T^{\frac{1}{1-m}}}u(\tau_1,0)
\ge \left[\frac{k_0\, M_0 }{2 T\, |B_{R_1}\setminus B_{2R_0}|}\right]^{\frac{1}{m}}\frac{t^{\frac{1}{1-m}}}{T^{\frac{1}{1-m}}}\\
&=\left[\frac{k_{s,d}}{4(R_1+R_0)^{d-2s}}\frac{(R_1-2R_0)^d}{R_1^d-(2R_0)^d}\right]^{\frac{1}{m}}
\frac{t^{\frac{1}{1-m}}}{T^{\frac{1}{m(1-m)}}}M_0^{\frac{1}{m}}\\
\end{split}
\end{equation}
since $t_*\le \tau_1\le T$\,. Moreover we have that
\begin{equation}\label{222}\begin{split}
u(t,0)
&\ge\left[\frac{k_{s,d}}{4(R_1+R_0)^{d-2s}}\frac{(R_1-2R_0)^{d-1}}{d(2R_0)^{d-1}}\right]^{\frac{1}{m}}
\frac{t^{\frac{1}{1-m}}}{T^{\frac{1}{m(1-m)}}}M_0^{\frac{1}{m}}\\
&=\left(\frac{k_{s,d}}{4d}\right)^{1/m}\left(\frac{R_1}{2R_0}-1\right)^{\frac{d-1}{m}}
\frac{t^{\frac{1}{1-m}}}{T^{\frac{1}{m(1-m)}}}\frac{M_0^{\frac{1}{m}}}{(R_1+R_0)^{\frac{d-2s}{m}}}\\
\end{split}
\end{equation}
where we have used the numerical inequality $a^d-b^d\le da^{d-1}(a-b)$\,, valid for any $a=R_1>2R_0=b$ to pass from \eqref{111} to \eqref{222}. By a standard argument it is easy to pass from the center to the infimum on $B_{R_0/2}(0)$ in the above estimates. The proof is concluded once we let $R_1=3R_0$.\qed

\medskip

\noindent {\bf Remarks. }(i) This result can be written alternatively as saying that there exists a universal constant $K_1=\max\{K^{-m}, C_*^{1/(1-m)}\}$ such for all solutions in the above class we have: for any $0\le t\le T$ and $R>0$
\begin{equation}\label{AC.s.VFDE}
\frac{\|u_0\|_{\LL^1(B_{R})}}{R^d }\le K_1\left[\frac{t^{\frac{1}{1-m}}}{R^{\frac{2s}{1-m}}} +\frac{T^{\frac{1}{1-m}}}
{t^{\frac{m}{1-m}}R^{2s}} \inf_{x\in B_{R/2}}u^m(t,x)  \right].
\end{equation}
This is easy to prove: by the previous Theorem, we have that either $t_*\le t$, that is
\[
\frac{\|u_0\|_{\LL^1(B_{R})}}{R^d}\le \left[\frac{t}{C_*R^{2s}}\right]^{\frac{1}{1-m}}
\]
or that $0\le t\le t_*$ and \eqref{thm.pos.T.VFDE} holds, namely
\[
\frac{\|u_0\|_{\LL^1(B_{R})}}{R^d } \le  \frac{T^{\frac{1}{1-m}}}
{K^{m} t^{\frac{m}{1-m}}R^{2s}} \inf_{x\in B_{R/2}}u^m(t,x)
\]
therefore, letting $K_1=\max\{K^{-m}, C_*^{1/(1-m)}\}$ we get \eqref{AC.s.VFDE}\,.

 This equivalent version is in complete formal agreement with  similar estimates proved by the authors in \cite{BV-ADV}, in the case $s=1$. However, our proof below  differs from the one in \cite{BV-ADV}, and provides an alternative proof when $s=1$. On the other hand, here we are considering solutions to the Cauchy problem, while in \cite{BV-ADV} we consider local weak solutions (i.e. without specifying boundary conditions). These estimates have been called Aronson-Caffarelli estimates in \cite{BV-ADV}, when $s=1$, since they are quite similar to the one that can be obtained for $m>1$, see Section \ref{sect.PME}. Finally we shall remark that in Section \ref{sect.ext.T} we will obtain quantitative upper estimates on the extinction time, and this will help to eliminate $T$ from the above lower estimates.

\noindent (ii) By comparison it is easy to prove that this estimates hold for a larger class of solutions, more precisely for the class of very weak solutions to the Cauchy Problem \eqref{FDE.eq} constructed in Theorem \ref{exist.large}, Section \ref{sect.exist.large}. This implies that the positivity result holds for solutions $u(t,\cdot)\in\LL^1(\RR^d, \varphi\dx)$ corresponding to initial data $0\le u_0\in\LL^1(\RR^d, \varphi\dx)$, where $\varphi$ is as in Theorem $\ref{prop.HP.s}$ with decay at infinity $|x|^{-\alpha}$, $d-[2s/(1-m)]<\alpha<d+(2s/m)$.

Once comparison is used, we can use as $T$ the extinction time of the reduced problem \ref{FDE.eq.Red.VFDE} in Step 1 of the above proof. In this way the quantitative result applies to solutions $u$ that may not extinguish in finite time.
Therefore we can interpret $T$ as the \textit{minimal life time for the solution $u(t,\cdot)$}, a concept that was already introduced by the authors in \cite{BV-ADV}\,, for which formula \eqref{t*.subcrit} provides a quantitative lower bound, namely
\begin{equation}\label{low.not.extinct}
t_*:= C_*\,R_0^{2s-d(1-m)}\|u_0\|_{\LL^1(B_{R_0})}^{1-m} \le T(u_0)\,.
\end{equation}
\begin{cor}[Solutions that do not extinguish in finite time]\label{cor.not.ext.}
Let $0<m<m_c$ and consider an initial datum $0\le u_0\in\LL^1(\RR^d, \varphi\dx)$, where $\varphi$ is as in Theorem $\ref{prop.HP.s}$, in particular, when $u_0\in\LL^1(\RR^d)$. Assume moreover that
\begin{equation}\label{555}
\liminf_{R\to +\infty}\;R^{\frac{2s}{1-m}-d}\|u_0\|_{\LL^1(B_R)}=+\infty\,.
\end{equation}
Then the corresponding solution $u(t,x)$ exists and is positive globally in space and time, hence does not extinguish in finite time. Moreover the quantitative lower bounds \eqref{thm.pos.T.VFDE} of Theorem $\ref{thm.lower.subcrit.}$ hold for any $0\le t\le t_*$ with $t_*$ given in \eqref{t*.subcrit} and $T=T(u_0\chi_{B_{R_0}})<+\infty$ is the extinction time of a reduced problem.
\end{cor}
\noindent {\sl Proof.~}If we consider an initial data with that behaviour at infinity, then by Theorem \ref{exist.large} there exists a very weak solution. By letting $R\to +\infty$ in the above lower bound \eqref{low.not.extinct} for $T$, to conclude that the minimal life time $T(u_0\chi_{B_{R}})\to \infty$, recalling that in this very fast diffusion range we have $2s<d(1-m)$, since $0<m<m_c$. \qed

\noindent\textbf{Remark. }A practical assumption on the initial datum $u_0$ that implies \eqref{555} is
\begin{equation}\label{555.1}
\liminf_{|x|\to +\infty}\;|x|^{\frac{2s}{1-m}}u_0(x)=+\infty\,.
\end{equation}
In view of Proposition $\ref{prop.ext}$ below, the exponent is sharp.

\subsection{Estimating the extinction time. }\label{sect.ext.T}

We next estimate the extinction time in terms of the initial data, extending a classical result of Benilan and Crandall \cite{BCr-cont}.  This is needed to eliminate the dependence on $T$ in the above lower estimates  when we consider initial data in $\LL^1(\RR^d)\cap\LL^{p_c}(\RR^d)$. For a detailed study of extinction time in the standard fast diffusion equation, see \cite{VazLN}.

\begin{prop}[Upper bounds for the extinction time]\label{prop.ext} Let $u$ be a weak solution to the equation \eqref{FDE.eq}, corresponding to $u_0\in \LL^1(\RR^d)\cap\LL^{p_c}(\RR^d)$ with $0<m<m_c=d/(d-2s)$, $0<s<1$ and let $p_c=d(1-m)/(2s)$. Then for all $0\le \tau\le t$ the following estimate holds true
\begin{equation}\label{ineq.pc}
\left[\int_{\RR^d}|u(t,x)|^{p_c}\dx\right]^{\frac{2s}{d}}\le \left[\int_{\RR^d}|u(\tau,x)|^{p_c}\dx\right]^{\frac{2s}{d}} - \frac{4m[d(1-m)-2s]}{d(d-2s)\mathcal{S}_s^2}(t-\tau)
\end{equation}
Moreover, there exists a finite extinction time $T\ge 0$ which can be bounded above as follows
\begin{equation}\label{ineq.T.ext}
T\le \frac{d(d-2s)\mathcal{S}_s^2}{4m[d(1-m)-2s]} \|u_0\|_{\LL^{p_c}(\RR^d)}^{1-m}\,.
\end{equation}
\end{prop}
\noindent {\sl Proof.~}The proof presented below is analogous to the one of Theorem 9.5 of \cite{DPQRV2}\,, but here we pay attention to the quantitative estimates\,. We multiply the equation by $|u|^{p-2}u$ with $p>1$, and integrate in $\RR^d$. Using
Strook-Varopoulos inequality \eqref{StrVar.ineq} in the form \eqref{StrVar.ineq.um}, we get
\begin{equation}\label{T.1.ext}\begin{split}
\frac{\rd}{\dt}\int_{\RR^d}|u(t,x)|^p\dx
    &=-p\int_{\RR^d}|u|^{p-2}u\,(-\Delta)^s (|u|^{m-1}u)\,\dx\\
    &\le -\frac{4mp(p-1)}{(p+m-1)^2}\int_{\RR^d}\left|(-\Delta)^{\frac{s}{2}}|u|^{\frac{p+m-1}{2}}\right|^2\dx\\
    &\le -\frac{4mp(p-1)}{(p+m-1)^2\mathcal{S}_s^2}\left(\int_{\RR^d}|u|^{\frac{d(p+m-1)}{d-2s}}\dx\right)^{\frac{d-2s}{d}}
\end{split}
\end{equation}
where in the last step we have used the Sobolev inequality \eqref{Sob.Fract} applied to $f=|u|^{\frac{p+m-1}{2}}$\,. Now we make the choice $p=p_c=d(1-m)/2s$, so that $p_c=\frac{d(p_c+m-1)}{d-2s}$, and we know that $p_c>1$ if and only if $m<m_c=d/(d-2s)$\,, and inequality \eqref{T.1.ext} becomes
\begin{equation}\label{T.2.ext}\begin{split}
\frac{\rd}{\dt}\int_{\RR^d}|u(t,x)|^{p_c}\dx
    &\le -\frac{4m[d(1-m)-2s]}{2s(d-2s)\mathcal{S}_s^2}\left(\int_{\RR^d}|u(t,x)|^{p_c}\dx\right)^{1-\frac{2s}{d}}
\end{split}
\end{equation}
Integrating the above differential inequality on $(\tau,t)$ gives both \eqref{ineq.pc} and inequality \eqref{ineq.T.ext}\,.\qed

Thanks to the above estimates we can get rid of the extinction time $T$ in the lower estimates of Theorem \ref{thm.lower.subcrit.}.
\begin{thm}[Local lower bounds II]\label{thm.lower.subcrit.Tpc}
Let $u$ be a weak solution to the equation \eqref{FDE.eq}, corresponding to $u_0\in \LL^1(\RR^d)\cap\LL^{p_c}(\RR^d)$ with $0<m<m_c=d/(d-2s)$, $0<s<1$ and let $p_c=d(1-m)/(2s)$. Then for every ball $B_{2R_0}\subset \Omega$, there exists a time
\begin{equation}\label{t*.subcrit.2}
t_*:= C_*\,R_0^{2s-d(1-m)}\|u_0\|_{\LL^1(B_{R_0})}^{1-m}\le T(u_0)\le \overline{C}\|u_0\|_{\LL^{p_c}(\RR^d)}^{1-m}\,,
\end{equation}
where we recall that $T(u_0)$ is the finite extinction time, such that
\begin{equation}\label{thm.pos.T.pc.VFDE}
\inf_{x\in B_{R_0/2}}u(t,x)\ge K_2 \frac{\|u_0\|_{\LL^1(B_{R_0})}^{\frac{1}{m}} }{R_0^{\frac{d-2s}{m}} } \frac{t^{\frac{1}{1-m}}}{\|u_0\|_{\LL^{p_c}(B_{R_0})}^{\frac{d}{2m}}}
\quad \mbox{ if } \ 0\le t\le t_*\,,
\end{equation}
where $C_*$ and $K_1$ are explicit positive universal constants, that depend only on $m,s,d$.
\end{thm}
The expression of the constants
\begin{equation}\label{const.subcrit.T.pc}\begin{split}
K_2 :=K\left[\frac{4m[d(1-m)-2s]}{d(d-2s)\mathcal{S}_s^2}\right]^{\frac{1}{m(1-m)}}\,,\qquad \overline{C}:=\frac{d(d-2s)\mathcal{S}_s^2}{4m[d(1-m)-2s]} \\
\end{split}\end{equation}\normalcolor
where $C_*$ and $K$ are as in \eqref{thm.pos.T.VFDE} and $k_{s,d}$ is the constant of the representation formula $  \varphi(x)=k_{s,d}\int_{\RR^d}\frac{\rho(y)}{|x-y|^{d-2s}}\dy$.

\noindent {\bf Remark.} This result can be written alternatively as saying that there exists a universal constant $K_3=\max\{K_2^{-m}, C_*^{1/(1-m)}\}$ such for all solutions in the above class we have: for any $0\le t\le T$ and $R>0$
\begin{equation}\label{AC.s.VFDE.pc}
\frac{\|u_0\|_{\LL^1(B_{R})}}{R^d \left[1\vee\|u_0\|_{\LL^{p_c}(B_{R})}\right]}
\le K_3\left[\frac{t^{\frac{1}{1-m}}}{R^{\frac{2s}{1-m}}} +\frac{1}{t^{\frac{m}{1-m}}R^{2s}} \inf_{x\in B_{R/2}}u^m(t,x)  \right]
\end{equation}
This equivalent version is in complete formal agreement with  similar estimates proved by the authors in \cite{BV-ADV}, in the case $s=1$.

\noindent {\sl Proof of Theorem \ref{thm.lower.subcrit.Tpc}.~}The proof is the same as for Theorem \ref{thm.lower.subcrit.}, as far as the first $3$ steps are concerned. At the end of Step 3, we need to bound from above the extinction time for the reduced problem \eqref{FDE.eq.Red.VFDE} with the estimates \eqref{ineq.T.ext} which give
\begin{equation}\label{ineq.T.ext.2}
T\le \frac{d(d-2s)\mathcal{S}_s^2}{4m[d(1-m)-2s]} \left[\int_{\RR^d}|u_0(x)|^{p_c}\dx\right]^{\frac{2s}{d}}=\frac{d(d-2s)\mathcal{S}_s^2}{4m[d(1-m)-2s]} \left[\int_{B_{R_0}}|u_0(x)|^{p_c}\dx\right]^{\frac{2s}{d}}\,,
\end{equation}
since $\supp(u_0)\subseteq B_{R_0}$\,. Then the proof follows simply by replacing $T$ with the above upper bound.\qed


\section{The Porous medium case}\label{sect.PME}

Lower estimates for nonnegative solutions of the standard porous medium equation were obtained in Aronson-Caffarelli in a famous paper \cite{ArCaff}. We want to show in this section how such a priori estimates extend to the fractional version considered in this paper.

\begin{thm}[Local lower bound]\label{thm.lower.pme}
 Let $u$ be a weak solution to Equation \eqref{FDE.eq}, corresponding to $u_0\in \LL^1(\RR^d)$. and let $m>1$.  We put $\vartheta:=1/[2s+d(m-1)]>0$. Then there exists a time
\begin{equation}\label{t*.PME}
t_*:=C \,R^{2s+d(m-1)}\,\|u_0\|_{\LL^1(B_{R})}^{-(m-1)}
\end{equation}
such that for every $t\ge t_*$ we have the lower bound
\begin{equation}
\inf_{x\in B_{R/2}}u(t,x)\ge K\,\dfrac{\|u_0\|_{\LL^1(B_{R})}^{2s\vartheta}}{t^{d\vartheta}}
\end{equation}
 valid for all  $R>0$. The positive constants $C$ and $K$  depend only on $m,s$ and $d$, and
 not on $R$.
\end{thm}

\noindent {\bf Remark.} This result can be written alternatively as saying that there exists a
universal constant $C_1=C_1(d,s,m)$ such for all solutions in the above class we have
\begin{equation}\label{AC.s.}
\int_{B_R(0)} u_0(x)\,dx \le C_1\left (R^{1/\vartheta(m-1)}t^{-1/(m-1)}+ u(0,t)^{1/2s\vartheta}t^{d/2s}\right).
\end{equation}
This equivalent version is in complete formal agreement with  Aronson-Caffarelli's estimate
for $s=1$. However, our proof below  differs very strongly from the ideas used in Aronson-Caffarelli's case since we cannot use the property of finite propagation of solution with compact support, which is false for $s<1$.

\noindent {\sl Proof.~} It is divided into several steps as follows.

\noindent$\bullet~$\textsc{Step 1. }\textit{Reduction. }By the comparison principle that it is sufficient to prove lower bounds for solutions $u$ to the following reduced problem:
\begin{equation}\label{FFDE.Prob.Red.PME}
\left\{
\begin{array}{lll}
\partial_t u
+(-\Delta)^s (u^m)=0\,,\; &\mbox{in }(0,\infty)\times\RR^d\,,\\
u(0,\cdot)=u_0\chi_{B_{R_0}}=\underline{u_0}\,,\; &\mbox{in }\RR^d\,,
\end{array}
\right.
\end{equation}
where $m>1$\,, $0<s<1$\,, and $R_0>0$\,. We only assume that $0\le u_0\in\LL^1(B_{R_0})$\,, which implies that $\underline{u_0}\in \LL^1(\RR^d)$ since $\supp(\underline{u_0})\subseteq B_{R_0}$ and also that $\|\underline{u_0}\|_{\LL^1(\RR^d)}=\|u_0\|_{\LL^1(B_{R_0})}$\,. It is not restrictive to assume that the ball $B_{R_0}$ is centered at the origin.

\noindent$\bullet~$\textsc{Step 2. }\textit{Smoothing effects. } In \cite{DPQRV2} there are the global $\LL^1$-$\LL^\infty$ smoothing effects,  which can be applied to solutions to our reduced Problem \ref{FFDE.Prob.Red.PME} as follows:\normalcolor
\begin{equation}\label{Smoothing.FPME}
\|u(t)\|_{\LL^\infty(\RR^d)}\le \frac{I_\infty}{t^{d\vartheta}}\|\underline{u_0}\|_{\LL^1(\RR^d)}^{2s\vartheta}
=\frac{I_\infty}{t^{d\vartheta}}\|u_0\|_{\LL^1(B_{R_0})}^{2s\vartheta}
\end{equation}
where $\vartheta=1/[2s+d(m-1)]$ and the constant $I_\infty$ only depends on $d,s,m$\,.

\noindent$\bullet~$\textsc{Step 3. }\textit{Aleksandrov principle. } We recall Theorem 11.2 of \cite{Vaz2012}. In view of the fact that the initial function is supported in the ball $B_{R_0}(0)$,  we have that
\[
u(t,0)\ge u(t,x)\,,\qquad\mbox{for all}\; t>0\;\mbox{and}\; |x|\ge 2R_0\,.
\]
Therefore, one has
\begin{equation}\label{aleks.1.PME.1}
\sup_{x\in \RR^d\setminus B_{2R_0}}u(t,x)\le u(t,0)\,.
\end{equation}

\noindent$\bullet~$\textsc{Step 4. }\textit{Weighted estimates. }If $\psi$ is a smooth, nonnegative, and sufficiently decaying function, we have
\[
\begin{split}
\left|\frac{\rd}{\dt}\int_{\RR^d}u(t,x)\psi(x)\dx\right|
&=\left|\int_{\RR^d}\left((-\Delta)^s u^m \right)\psi\dx\right|
=_{(a)}\left|\int_{\RR^d} u^m (-\Delta)^s\psi\dx\right|\\
&\le \|u(t)\|_{\LL^\infty(\RR^d)}^{m-1}\left\|(-\Delta)^s\psi\right\|_{\LL^\infty(\RR^d)} \int_{\RR^d} u(t,x) \dx\\
&\le_{(b)} \frac{I_\infty^{m-1}}{t^{d\vartheta(m-1)}}\|u_0\|_{\LL^1(B_{R_0})}^{2s\vartheta(m-1)}
    \left\|(-\Delta)^s\psi\right\|_{\LL^\infty(\RR^d)} \int_{\RR^d} \underline{u_0}(x) \dx\\
&:= \frac{I_\infty^{m-1}}{t^{d\vartheta(m-1)}}
    \left\|(-\Delta)^s\psi\right\|_{\LL^\infty(\RR^d)}\|u_0\|_{\LL^1(B_{R_0})}^{2s\vartheta(m-1)+1}
:= \frac{K[u_0,\psi]}{t^{d\vartheta(m-1)}}\,.
\end{split}
\]
Notice that in $(a)$ we have used the fact that $(-\Delta)^s$ is a symmetric operator, while in $(b)$ we have used the smoothing effect \eqref{Smoothing.FPME} of Step 2 and the the conservation of mass:  $\int_{\RR^d} u(t,x) \dx=\int_{\RR^d} \underline{u_0}(x) \dx$\,, for all $t>0$, together with the fact that $\supp(\underline{u_0})\subseteq B_{R_0}$.
We refer to \cite{DPQRV2} for a proof of the smoothing effect and of the conservation of mass.  Summing up,
\[
\left|\frac{\rd}{\dt}\int_{\RR^d}u(t,x)\psi(x)\dx\right|\le \frac{K[u_0,\psi]}{t^{1-2s\vartheta}}\,,
\]
since $d\vartheta(m-1)=1-2s\vartheta$. Integrating the above differential inequality on $(0,t)$ with $t\ge 0$ we obtain:
\[
-\frac{K[u_0,\psi]}{2s\vartheta}\,t^{2s\vartheta}\le \int_{\RR^d}u(t,x)\psi(x)\dx-\int_{\RR^d}u(0,x)\psi(x)\dx
\le \frac{K[u_0,\psi]}{2s\vartheta}\,t^{2s\vartheta}\,.
\]
We will use this in the form
\begin{equation}\label{step4.1}
\int_{\RR^d}u(0,x)\psi(x)\dx - \frac{K[u_0,\psi]}{2s\vartheta}\,t^{2s\vartheta}\le \int_{\RR^d}u(t,x)\psi(x)\dx\,.
\end{equation}
Moreover, if $\psi\in\LL^1(\RR^d)$ and $R_1\ge 2R_0$, we have
\begin{equation}\label{step4.3}\begin{split}
\int_{\RR^d}u(t,x)\psi(x)\dx
    &= \int_{B_{R_1}}u(t,x)\psi(x)\dx+ \int_{B_{R_1}^c}u(t,x)\psi(x)\dx\\
    &\le_{(a)} |B_{R_1}|\sup_{|x|\le {R_1}}u(t,x) + \sup_{x\in \RR^d\setminus B_{2R_0}}u(t,x)\int_{B_R^c}\psi(x)\dx\\
        &\le_{(b)} |B_{R_1}|\frac{I_\infty}{t^{d\vartheta}}\|u_0\|_{\LL^1(B_{R_0})}^{2s\vartheta} + u(t,0)\int_{B_{R_1}^c}\psi(x)\dx
\end{split}
\end{equation}
where in $(a)$ we have used the fact that $\psi\le 1$, $R_1\ge 2R_0$\,, and that $\psi\in\LL^1(\RR^d)$. In $(b)$ we have used the smoothing effect \eqref{Smoothing.FPME} of step 2 and the Aleksandrov principle of Step 3\,.
Putting together inequalities \eqref{step4.1} and \eqref{step4.3}\,, we obtain
\begin{equation}\label{step4.4}
\int_{\RR^d}u(0,x)\psi(x)\dx - \frac{K[u_0,\psi]}{2s\vartheta}\,t^{2s\vartheta}-|B_{R_1}|\frac{I_\infty}{t^{d\vartheta}}\|u_0\|_{\LL^1(B_{R_0})}^{2s\vartheta} \le u(t,0)\int_{B_{R_1}^c}\psi(x)\dx\,.
\end{equation}
Next, in order to estimate $K[u_0,\psi]$ in a convenient way we take $\psi(x)=\phi(|x|/R)$ with $\phi$ as in Lemma \ref{Lem.phi}\,, we have $|(-\Delta)^s\psi|\le c_3 R^{-2s}$\,, for some  constant $c_3= c_3(d,s)$\,. Then,
\begin{equation}\label{step4.2}\begin{split}
K[u_0,\psi]=I_\infty^{m-1}\left\|(-\Delta)^s\psi\right\|_{\LL^\infty(\RR^d)}\|u_0\|_{\LL^1(B_{R_0})}^{2s\vartheta(m-1)+1}
\le \frac{c_3 I_\infty^{m-1}}{R^{2s}}\|u_0\|_{\LL^1(B_{R_0})}^{2s\vartheta(m-1)+1}
\end{split}\end{equation}
 When $R\ge R_0$ and $R_1\ge 2R_0$, we arrive at
\begin{equation}\label{step4.5}\begin{split}
\|u_0\|_{\LL^1(B_{R_0})} - \frac{c_3 I_\infty^{m-1}}{2s\vartheta}\frac{\|u_0\|_{\LL^1(B_{R_0})}^{2s\vartheta(m-1)+1}}{R^{2s}}
\,t^{2s\vartheta}
&-\omega_d R_1^d \frac{I_\infty}{t^{d\vartheta}}\|u_0\|_{\LL^1(B_{R_0})}^{2s\vartheta} \le u(t,0)\int_{B_{R_1}^c}\psi(x)\dx\\
&\le u(t,0)R^d \int_{\RR^d}\varphi(x)\dx=c_4 R^d u(t,0)\,.
\end{split}
\end{equation}

\noindent$\bullet~$\textsc{Step 5. }\textit{Choosing the parameters. }We want to choose $t>0$, $R\ge R_0$ and $R\ge 2R_0$ so that the left-hand side of \eqref{step4.5} is larger than $\|u_0\|_{\LL^1(B_{R_0})}/2$, which will then give the desired bound from below for $u(0,t)$.
We first make the choice
\begin{equation}
R_1^d=\frac{\|u_0\|_{\LL^1(B_{R_0})}^{d(m-1)\vartheta}}{4 \omega_d I_\infty}{t^{d\vartheta}}\,,
\end{equation}
which will satisfy the condition $R_1\ge 2R_0$ if and only if $t\ge t_*$ where
\begin{equation}
t_*^{\vartheta}=c_5 R_0\|u_0\|_{\LL^1(B_{R_0})}^{-(m-1)\vartheta}, \qquad c_5=2^{1+(2/d}(\omega_d I_\infty)^{1/d}.
\end{equation}
Now we can make the second choice, $R$ has to be large enough, for instance:
\begin{equation}
R=c_6\|u_0\|_{\LL^1(B_{R_0})}^{\vartheta(m-1)}\,t^{\vartheta}\,,\qquad c_6=\left(\frac{4\,c_3 I_\infty^{m-1}}{2s\vartheta}\right)^{1/2s}
\end{equation}
Both choices will give for $t\ge t_*$ the lower bound
\[
\frac{\|u_0\|_{\LL^1(B_{R_0})}}{2c_4\, R^d}\le  u(t,0)\,,
\]
which can be rewritten as
\begin{equation}\label{step4.6}
c_7\frac{\|u_0\|_{\LL^1(B_{R_0})}^{2s\vartheta}}{t^{d\vartheta}}
\le  u(t,0)\,,\quad \mbox{for any}\quad t\ge t_*\,,\quad \mbox{with}\quad c_7=\frac{1}{2c_4\, c_6^d}\mbox{\,.}
\end{equation}
By a standard argument it is easy to pass from the center to the infimum on $B_{R_0/2}(0)$ in the above estimates.\qed

\medskip

\noindent {\bf Remark}. In the limit $m\to 1$ of the estimate of Theorem \ref{thm.lower.pme} we  obtain the result  of  Proposition \ref{thm.lower.m=1} for $m=1$.

\medskip

 \noindent {\bf Open Problem.} To calculate the positivity of the solutions for small times is not known yet.


\section{Existence and uniqueness of initial traces}\label{sect.traces}

The existence of solutions of the Cauchy Problem (\ref{FDE.eq})-(\ref{FDE.id}) can be extended to the case where the initial datum is a finite and nonnegative Radon measure. We denote by
 ${\cal M}^+(\ren)$ the space of such measures on $\ren$. Here is the result proved in Theorem 4.1 of \cite{Vaz2012}.

\noindent {\bf Theorem.} {\sl For every $\mu\in {\cal M}^+(\ren)$ there exists a nonnegative and continuous weak solution of Equation \eqref{FDE.eq} in $Q=(0,\infty)\times \RR^d$ taking initial data $\mu$ in the sense that for every $\varphi\in C_c^{2}(\RR^d)$ we have
\begin{equation}
\lim_{t\to 0^+} \int u(t,x)\varphi(x)\,dx=\int \varphi(x)d\mu(x)\,.
\end{equation}
}
In this section we address the reverse problem, i.\,e., given a solution to find the initial trace. In the case $s=1$ such question was solved thanks to the works of Aronson-Caffarelli \cite{ArCaff}, Dahlberg-Kenig \cite{DK2}, Pierre \cite{Pierre} and others, see a presentation in \cite{VazBook}, Chapter 13.

\begin{lem}[Conditions for existence and uniqueness of initial traces]\label{lem.init.trace}
Let $m>0$ and let $u$ be a solution to  equation \eqref{FDE.eq} in $(0,T]\times\RR^d$. Assume that there exist a time $0<T_1\le T$, some positive constants $K_1,K_2,\alpha>0$  and a continuous function $\omega:[0,+\infty)\to [0,+\infty)$, with $\omega(0)=0$ such that
\begin{equation}\label{hyp.1.lem.init}
(i)\qquad \sup_{t\in (0,T_1]}\int_{B_R(x_0)}u(t,x)\dx\le K_1\,,\qquad\forall\; R>0\,,\;x_0\in\RR^d\,,
\end{equation}
as well as
\begin{equation}\label{hyp.2.lem.init}
(ii)\qquad \left[\int_{\RR^d}u(t,x)\varphi(x)\dx\right]^\alpha\le \left[\int_{\RR^d}u(t',x)\varphi(x)\dx\right]^\alpha
    +K_2\,\omega(|t-t'|)
\end{equation}
for all $0<t,t'\le T_1$ and for all $\varphi\in C_c^\infty(\RR^d)$\,. Then there exists a unique nonnegative Radon measure $\mu$ as initial trace, that is
\[
\int_{\RR^d}\varphi\,\rd\mu=\lim_{t\to 0^+}\int_{\RR^d}u(t,x)\varphi(x)\dx\,,\qquad\mbox{for all }\varphi\in C_c^\infty(\RR^d)\,.
\]
Moreover the initial trace $\mu$ satisfies the bound \eqref{hyp.1.lem.init} with the same constant, namely $\mu(B_{R}(x_0))\le K_1$\,. \end{lem}

Notice that the constants $K_1$ and $K_2$ may depend on $u$ and $\varphi$, usually through some norm.

\medskip

\noindent {\sl Proof.~}The proof is divided in two steps in which we prove existence and uniqueness of the initial trace respectively.

\noindent$\bullet~$\textsc{Step 1. }\textit{Existence of the initial trace. }Hypothesis $(i)$ easily implies that
\[
\limsup_{t\to 0^+}\int_{B_R(x_0)}u(t,x)\dx\le K_1\,,\qquad\forall\; R>0\,,\;x_0\in\RR^d\,.
\]
Moreover, it implies weak compactness for measures (to be more precise, weak$^*$ compactness, see Theorem \ref{thm.cpt.radon} in the Appendix \ref{appendix.meas}), so that there exists a sequence $t_k\to 0^+$ as $k\to \infty$ with $0<t_k<T_1$\,, and a nonnegative Radon measure $\mu$ so that
\[
\lim_{k\to\infty}\int_{\RR^d}u(t_k,x)\varphi(x)\,\dx=\int_{\RR^d}\varphi\,\rd\mu\qquad\mbox{for all }\varphi\in C^0_c(\RR^d)\,.
\]
The bound \eqref{intit.trace.bdd.lem} on the initial trace: $\mu(B_R(x_0))\le K_1$ follows from the above bound on the $\limsup$\,.

\noindent$\bullet~$\textsc{Step 2. }\textit{Uniqueness of the initial trace. }  The initial trace whose existence we have just proved may, of course, depend on the sequence $t_k$. We will now show  that this is not the case, thanks to hypothesis $(ii)$. Assume that there exist two sequences $t_k\to 0^+$ and $t'_k\to 0^+$ as $k\to \infty$\,, so that $u(t_k)\to \mu$ and $u(t'_k)\to \nu$,with $\mu, \nu\in {\cal M}^+(\ren)$. We will prove that
\begin{equation}\label{step.2.uniq.1}
\int_{\RR^d}\varphi\,\rd\mu=\int_{\RR^d}\varphi\,\rd\nu\qquad\mbox{for all }\varphi\in C^{\infty}_c(\RR^d)\,.
\end{equation}
so that $\mu=\nu$ as positive linear functionals on $C^{\infty}_c(\RR^d)$. Then by the Riesz representation theorem (cf. Theorem \ref{thm.riesz.radon}) we know that $\mu=\nu$ also as Radon measures on $\RR^d$. Therefore, it only remains to prove \eqref{step.2.uniq.1}\,: hypothesis $(ii)$ implies that for any $t,t' >0$\,, with $0<t+t'\le T_1\le T$, and any $\varphi\in C^{\infty}_c(\RR^d)$\, we have $\omega(|(t+t')-t|)=\omega(t')$ and
\begin{equation}\label{lem61.step2.1}
\left[\int_{\RR^d}u(t,x)\varphi(x)\dx\right]^\alpha\le \left[\int_{\RR^d}u(t+t',x)\varphi(x)\dx\right]^\alpha
    +K_2\omega(t')\,.
\end{equation}
First we let $t=t_k$ and $t'>0$ to be chosen later, then we let $t_k\to 0^+$ so that $u(t_k)\rightharpoonup \mu$, and we get
\begin{equation}\label{lem61.step2.2}
\left[\int_{\RR^d}\varphi\,\rd\mu\right]^\alpha\le \left[\int_{\RR^d}u(t',x)\dx\right]^\alpha
    +K_2\omega(t')\,.
\end{equation}
Then we put $t'=t'_k$ and let $t'_k\to 0^+$ so that $u(t'_k)\rightharpoonup \nu$, $\omega(t'_k)\to \omega(0)=0$ and we obtain the first inequality
\begin{equation}\label{lem61.step2.3}
\left[\int_{\RR^d}\varphi\,\rd\mu\right]^\alpha\le\left[\int_{\RR^d}\varphi\,\rd\nu\right]^\alpha\,.
\end{equation}
Then, we proceed exactly in the same way but we exchange the roles of $t_k$ and $t'_k$ to obtain the opposite inequality $\left[\int_{\RR^d}\varphi\,\rd\mu\right]^\alpha\le\left[\int_{\RR^d}\varphi\,\rd\nu\right]^\alpha\,.$ Therefore we have that $\mu=\nu$ as positive linear functionals on $C^{\infty}_c(\RR^d)$\, as desired.\qed

\begin{thm}[Existence and uniqueness of initial trace, FD case]\label{thm.init.trace.m<1}
Let $0<m<1$ and let $u$ be a nonnegative weak solution of equation \eqref{FDE.eq} in $(0,T]\times\RR^d$. Assume that $\|u(T)\|_{\LL^1(\RR^d)}<\infty$. Then there exists a unique nonnegative Radon measure $\mu$ as initial trace, that is
\begin{equation}\label{eq.trace1}
\int_{\RR^d}\psi\,\rd\mu=\lim_{t\to 0^+}\int_{\RR^d}u(t,x)\psi(x)\,\dx\,,\qquad\mbox{for all }\psi\in C_0(\RR^d)\,.
\end{equation}
Moreover, the initial trace $\mu$ satisfies the bound
\begin{equation}\label{intit.trace.bdd.lem}
\mu(B_{R}(x_0))\le \|u(T)\|_{\LL^1(\RR^d)} + C_1 R^{d(1-m)-2s}\,T\,.
\end{equation}
where $C_1=C_1(m,d,s)>0$ as in \eqref{HP.s}.
\end{thm}
\noindent {\sl Proof.~}The proof is divided into three steps.

\noindent$\bullet~$\textsc{Step 1. }\textit{Weighted estimates I. Existence. }First we recall the weighted estimates of Theorem \ref{prop.HP.s}\,, which imply for all $0\le t\le T_1\le T$
\begin{equation}\label{HP.s.2.mu}\begin{split}
\left(\int_{\RR^d}u(t,x)\phi_R(x)\dx\right)^{1-m}
&\le\left(\int_{\RR^d}u(T,x)\phi_R(x)\dx\right)^{1-m}
+ C_1 R^{d(1-m)-2s}\,|T-T_1|\\
&\le  \|u(T)\|_{\LL^1(\RR^d)}  + C_1 R^{d(1-m)-2s}\,T := K_1
\end{split}
\end{equation}
since $\phi_R\le 1$ and where $C_1>0$ depends only on $\alpha,m,d$ as in Theorem \ref{prop.HP.s}\,. Since $\phi_R\ge 1$ on $B_R$ it is clear that this implies hypothesis $(i)$ of Lemma \ref{lem.init.trace}, therefore it guarantees the existence of an initial trace that satisfies the bound $\mu(B_R(x_0))\le K_1=\|u(T)\|_{\LL^1(\RR^d)} + C_1 R^{d(1-m)-2s}\,T$\,.

\noindent$\bullet~$\textsc{Step 2. }\textit{Pseudo-local estimates. Uniqueness. }In order to prove uniqueness of the initial trace is is sufficient to prove hypothesis $(ii)$ of Lemma \ref{lem.init.trace}, namely we need to prove that
\begin{equation}\label{hyp.2.lem.init.2}
\left[\int_{\RR^d}u(t,x)\psi(x)\dx\right]^\alpha\le \left[\int_{\RR^d}u(t',x)\psi(x)\dx\right]^\alpha
    +K_2\,\omega(|t-t'|)
\end{equation}
for all $0<t,t'\le T_1\le T$ and for all $\psi\in C_c^\infty(\RR^d)$\,. We will see that this is true for $\alpha=1$ and $\omega(|t-t'|)=|t-t'|$. Let $\psi\in C_c^\infty(\RR^d)$, then we have
\[
\begin{split}
\left|\frac{\rd}{\dt}\int_{\RR^d}u(t,x)\psi(x)\dx\right|
&=\left|\int_{\RR^d}(-\Delta)^s u^m\psi\dx\right|
=_{(a)}\left|\int_{\RR^d}u^m(-\Delta)^s\psi\dx\right|
\le \int_{\RR^d}u^m\,\phi_R(x)\,\frac{\left|(-\Delta)^s\psi(x)\right|}{\phi_R(x)}\dx\\
&\le_{(b)}  \left\|\frac{\left|(-\Delta)^s\psi(x)\right|}{\phi_R(x)}\right\|_{\LL^\infty(\RR^d)}
      \left(\int_{\RR^d}\phi_R\dx\right)^{1-m}\,\left(\int_{\RR^d}u\phi_R\dx\right)^m        \\
&\le k_7\,\|\phi_R\|_{\LL^1(\RR^d)}\,K_1:= K_2\,.
\end{split}
\]
Notice that in $(a)$ we have used the fact that $(-\Delta)^s$ is a symmetric operator. In $(b)$ we have chosen $\phi_R(x):=\phi(x/R)$\,, with $\phi$ as in \eqref{phi} of Lemma \ref{Lem.phi}, with the decay at infinity $\alpha=d+2s$. It then follows that
\[
\left\|\frac{\left|(-\Delta)^s\psi(x)\right|}{\phi_R(x)}\right\|_{\LL^\infty(\RR^d)}\le k_7\,,
\]
since we know by Lemma \ref{Lem.phi} that $\left|(-\Delta)^s\psi(x)\right|\le k_5\,|x|^{-(d+2s)}$\,, \normalcolor and we have chosen $\phi_R\ge k_6/|x|^{d+2s}$\,.
We have also used the fact that the $\LL^m$-norm ($m<1$) is less than the $\LL^1$ norm since the measure $\phi_R\dx$ is finite. In the last line of the display we have used the bound of Step 1, namely that $\left(\int_{\RR^d}u(t,x)\varphi_R(x)\dx\right)^{1-m}\le M_T + C_1 R^{d(1-m)-2s}\,T := K_1$ for all $0\le t\le T_1$\,. Summing up, we  have obtained:
\[
\left|\frac{\rd}{\dt}\int_{\RR^d}u(t,x)\psi(x)\dx\right|\le K_2\,.
\]
Integrating the above differential inequality we obtain:
\begin{equation}\label{est.loc.nonloc.0}
\int_{\RR^d}u(t,x)\psi(x)\dx\le \int_{\RR^d}u(\tau,x)\psi(x)\dx
+ K_2\,|t-\tau|\qquad\mbox{ for any $\tau,t\ge 0$ and all $\psi\in C_c^\infty(\RR^d)$\,.}
\end{equation}
\noindent$\bullet~$\textsc{Step 3. } We still have to pass from test functions $\psi\in C_c^\infty(\RR^d)$ to $\psi\in C_c^0(\RR^d)$  in formula \eqref{eq.trace1}, but this is easy by approximation (mollification).\qed

\noindent\textbf{Remarks. } (i) The proof applies with minor modification to the class of solutions with data $u_0\in \LL^1(\RR^d,\varphi\dx)$ constructed in Section \ref{sect.exist.large}\,.

\noindent(ii) Notice that estimates \eqref{est.loc.nonloc.0} are only pseudo-local estimates: the global information about $u(T)$, namely the bound $\|u(T)\|_{\LL^1(\RR^d)}$ is contained in the constant $K_1$ and therefore in $K_2$.

\noindent (iii) The existence of solutions and traces for the standard FDE with (not necessarily locally finite) Borel measures as data is studied in Chasseigne-Vazquez \cite{ChVaz02}. We do not address the corresponding question here.


\color{darkblue}

\begin{thm}[\color{darkblue}\bf Existence and uniqueness of initial trace, HE case $m=1$]\label{thm.init.trace.m=1}
Let $m=1$ and let $u$ be a nonnegative weak solution of equation \eqref{FDE.eq} in $(0,T]\times\RR^d$. Assume that $\|u(T)\|_{\LL^1(\RR^d,\varphi)}<\infty$ where $\varphi$ is as in Theorem $\ref{prop.HP.s}$ with decay at infinity $|x|^{-\alpha}$, $\alpha=d+2s$. Then there exists a unique nonnegative Radon measure $\mu$ as initial trace, that is
\begin{equation*}
\int_{\RR^d}\psi\,\rd\mu=\lim_{t\to 0^+}\int_{\RR^d}u(t,x)\psi(x)\,\dx\,,\qquad\mbox{for all }\psi\in C_0(\RR^d)\,.
\end{equation*}
Moreover, the initial trace $\mu$ satisfies the bound
\begin{equation*}
\mu(B_{R}(x_0))\le \ee^{K_0\,T}\,\|u(T)\|_{\LL^1(\RR^d,\varphi)}\,,
\end{equation*}
where $K_0=K_0(m,d,s)>0$ as in \eqref{HE.2}.
\end{thm}

\noindent{\bf Proof.~}The proof is divided into three steps.

\noindent$\bullet~$\textsc{Step 1. }\textit{Weighted estimates I. Existence. }First we prove the weighted estimates when $m=1$: for all $0\le t\le T_1\le T$
\begin{equation}\label{HE.1}\begin{split}
\int_{\RR^d}u(t,x)\phi_R(x)\dx
&\le\ee^{K_0(T-t)}\,\int_{\RR^d}u(T,x)\phi_R(x)\dx
\le \ee^{K_0\,T}\,\|u(T)\|_{\LL^1(\RR^d,\phi_R)}:= K_1
\end{split}
\end{equation}
The proof of the above inequality is as follows. Consider a function  $\varphi=\phi_R$ as in Theorem $\ref{prop.HP.s}$ with decay at infinity $|x|^{-\alpha}$, $\alpha=d+2s$ and such that $\phi_R\ge 1$ on $B_R$, so that by Lemma \eqref{Lem.phi} we have
\begin{equation}\label{HE.2}
\left\|\frac{(-\Delta)^s\phi_R}{\phi_R}\right\|_\infty \le K_0<+\infty
\end{equation}
so that
\[
\left|\frac{\rd}{\dt}\int_{\RR^d}u\phi_R \dx \right|
= \left|\int_{\RR^d}u\, (-\Delta)^s\phi_R\dx \right|
\le \left\|\frac{(-\Delta)^s\phi_R}{\phi_R}\right\|_\infty\,\int_{\RR^d} u\, \phi_R\dx\le K_0 \int_{\RR^d} u\, \phi_R\dx\,.
\]
from which \eqref{HE.1} follows. Since $\phi_R\ge 1$ on $B_R$ it is clear that \eqref{HE.1} implies hypothesis $(i)$ of Lemma \ref{lem.init.trace}, therefore it guarantees the existence of an initial trace that satisfies the bound
$\mu(B_R(x_0))\le K_1=\ee^{K_0\,T}\,\|u(T)\|_{\LL^1(\RR^d,\varphi)}\,.$

\noindent$\bullet~$\textsc{Step 2. }\textit{Pseudo-local estimates. Uniqueness. }In order to prove uniqueness of the initial trace is is sufficient to prove hypothesis $(ii)$ of  Lemma \ref{lem.init.trace}, namely we need to prove that
\begin{equation}\label{HE.3}
\left[\int_{\RR^d}u(t,x)\psi(x)\dx\right]^\alpha\le \left[\int_{\RR^d}u(t',x)\psi(x)\dx\right]^\alpha
    +K_2\,\omega(|t-t'|)
\end{equation}
for all $0<t,t'\le T_1\le T$ and for all $\psi\in C_c^\infty(\RR^d)$\,. We will see that this is true for $\alpha=1$ and $\omega(|t-t'|)=|t-t'|$.   Let $\psi\in C_c^\infty(\RR^d)$, then we have
\[
\begin{split}
\left|\frac{\rd}{\dt}\int_{\RR^d}u(t,x)\psi(x)\dx\right|
&=\left|\int_{\RR^d}(-\Delta)^s u\psi\dx\right|
=_{(a)}\left|\int_{\RR^d}u^m(-\Delta)^s\psi\dx\right|
\le \int_{\RR^d}u\,\phi_R(x)\,\frac{\left|(-\Delta)^s\psi(x)\right|}{\phi_R(x)}\dx   \\
&\le_{(b)}  \left\|\frac{\left|(-\Delta)^s\psi(x)\right|}{\phi_R(x)}\right\|_{\LL^\infty(\RR^d)}
       \int_{\RR^d}u\phi_R\dx
\le K'_0\,K_1:= K_2\,.
\end{split}
\]
Notice that in $(a)$ we have used the fact that $(-\Delta)^s$ is a symmetric operator.   In $(b)$ we have chosen $\phi_R(x)$\,, as in Step 1, so that, by Lemma \eqref{Lem.phi}, we have
\[
\left\|\frac{\left|(-\Delta)^s\psi(x)\right|}{\phi_R(x)}\right\|_{\LL^\infty(\RR^d)}\le K'_0\,.
\]
Finally, we have used the bound \eqref{HE.1} of Step 1: for all $0\le t\le T_1$ we have $\int_{\RR^d}u(t,x)\phi_R(x)\dx\le  K_1$ \,.
Summing up, we  have obtained that for all $t\ge 0$ and all $\psi\in C_c^\infty(\RR^d)$
\[
\left|\frac{\rd}{\dt}\int_{\RR^d}u(t,x)\psi(x)\dx\right|\le K_2\,.
\]
Integrating the above differential inequality, we obtain \eqref{HE.3} with $\alpha=1$ and $\omega(|t-t'|)=|t-t'|$\,.

\noindent$\bullet~$\textsc{Step 3. } We still have to pass from test functions $\psi\in C_c^\infty(\RR^d)$ to $\psi\in C_c^0(\RR^d)$  in formula \eqref{eq.trace1}, but this is easy by approximation (mollification).\qed


\normalcolor

\begin{thm}[Existence and uniqueness of initial trace, PME case]\label{thm.init.trace.m>1}
Let $m>1$ and let $u$ be a solution to the Cauchy problem \ref{FDE.eq} on $(0,T]\times\RR^d$. Assume that $\|u(T)\|_{\LL^1(\RR^d)}+\|u(T)\|_{\LL^\infty(\RR^d)}<+\infty$. Then there exists a unique nonnegative Borel measure $\mu$ as initial trace, that is
\begin{equation}\label{eq.trace2}
\int_{\RR^d}\psi\,\rd\mu=\lim_{t\to 0^+}\int_{\RR^d}u(t,x)\psi(x)\dx\,,\qquad\mbox{for all }\psi\in C_0(\RR^d)\,.
\end{equation}
Moreover the initial trace $\mu$ satisfies the bound
\begin{equation}\label{intit.trace.bdd.lem.1}
\mu(B_{R}(x_0))\le C_1\left[\left(\frac{R^{2s+d(m-1)}}{T}\right)^{\frac{1}{m-1}}+ T^{\frac{d}{2s}}\, u(x_0,T)^{\frac{1}{2s\vartheta}}\right]\,,
\end{equation}
where $C_1=C_1(m,d,s)>0$ as in Theorem \ref{thm.lower.pme}.
\end{thm}
\noindent {\sl Proof.~}The proof is divided in three steps

\noindent$\bullet~$\textsc{Step 1. }\textit{Weighted estimates I. Existence. }First we recall the lower bounds of Theorem \eqref{thm.lower.pme} rewritten in the form \eqref{AC.s.}
\begin{equation}\label{AC.s.2}\begin{split}
\int_{B_R(x_0)} u(\tau,x)\,dx
&\le C_1\left[\left(\frac{R^{2s+d(m-1)}}{T}\right)^{\frac{1}{m-1}}+ T^{\frac{d}{2s}}\, u(x_0,T)^{\frac{1}{2s\vartheta}}\right]
:= K_1.
\end{split}
\end{equation}
on the time interval $(\tau,T]\subseteq (0,T]$\,. It is clear that this implies hypothesis $(i)$ of Lemma \ref{lem.init.trace}, therefore it guarantees the existence of an initial trace that satisfy the bound $\mu(B_R(x_0))\le K_1\,.$

\noindent$\bullet~$\textsc{Step 2. }\textit{Smoothing effects and mass conservation. }In \cite{DPQRV2} there are the global $\LL^1-\LL^\infty$ smoothing effects which provide global upper bounds for solutions to the Cauchy problem \ref{FDE.eq}\,. We apply such smoothing effects to solutions to our reduced Problem \ref{FFDE.Prob.Red} to get
\begin{equation}\label{Smoothing.FPME.1}
\|u(t)\|_{\LL^\infty(\RR^d)}\le \frac{2^{d\vartheta}I_\infty}{t^{d\vartheta}}\|u(t/2)\|_{\LL^1(\RR^d)}^{2s\vartheta}
\end{equation}
where $\vartheta=1/[2s+d(m-1)]$ and the constant $I_\infty$ only depends on $d,s,m$\,. Moreover, we know that there holds also the conservation of mass on the time interval $[t/2,T]\subset (0,T]$, so that inequality \eqref{Smoothing.FPME.1} becomes
\begin{equation}\label{Smoothing.FPME.2}
\|u(t)\|_{\LL^\infty(\RR^d)}\le \frac{2^{d\vartheta}I_\infty}{t^{d\vartheta}}\|u(T)\|_{\LL^1(\RR^d)}^{2s\vartheta}\,.
\end{equation}

\noindent$\bullet~$\textsc{Step 3. }\textit{Weighted estimates II. Pseudo-local estimates. Uniqueness. }In order to prove uniqueness of the initial trace is is sufficient to prove hypothesis $(ii)$ of Lemma \ref{lem.init.trace}, namely we need to prove
\begin{equation}\label{hyp.2.lem.init.2.2}
\left[\int_{\RR^d}u(t,x)\psi(x)\dx\right]^\alpha\le \left[\int_{\RR^d}u(t',x)\psi(x)\dx\right]^\alpha
    +K_2\,\omega(|t-t'|)
\end{equation}
for all $0<t,t'\le T_1\le T$ and for all $\psi\in C_c^\infty(\RR^d)$\,. We will see that this is true for $\alpha=1$ and $\omega(|t-t'|)=|t^\sigma-t'^\sigma|)$ with $\sigma=2s/[2s+d(m-1)]$\,. Let $\psi\in C_c^\infty(\RR^d)$, then we have
\[
\begin{split}
\left|\frac{\rd}{\dt}\int_{\RR^d}u(t,x)\psi(x)\dx\right|
&=\left|\int_{\RR^d}(-\Delta)^s u^m\psi\dx\right|
=_{(a)}\left|\int_{\RR^d}u^m(-\Delta)^s\psi\dx\right|\\
&\le \|u(t)\|_{\LL^\infty(\RR^d)}^{m-1} \left\|(-\Delta)^s\psi\right\|_{\LL^\infty(\RR^d)} \int_{\RR^d}u(t) \dx\\
&\le_{(b)} \frac{2^{d\vartheta(m-1)}I_\infty^{m-1}}{t^{d\vartheta(m-1)}}\|u(T)\|_{\LL^1(\RR^d)}^{2s\vartheta(m-1)} \left\|(-\Delta)^s\psi\right\|_{\LL^\infty(\RR^d)}\|u(t)\|_{\LL^\infty(\RR^d)}\\
&=\frac{2^{d\vartheta(m-1)}I_\infty^{m-1}}{t^{d\vartheta(m-1)}}\|u(T)\|_{\LL^1(\RR^d)}^{2s\vartheta(m-1)+1} \left\|(-\Delta)^s\psi\right\|_{\LL^\infty(\RR^d)}
:= \frac{K_2}{t^{d\vartheta(m-1)}}\,.
\end{split}
\]
Notice that in $(a)$ we have used the fact that $(-\Delta)^s$ is a symmetric operator. In $(b)$ we have used the smoothing effect \eqref{Smoothing.FPME.2} of Step 2\,. Summing up we  have obtained:
\[
\left|\frac{\rd}{\dt}\int_{\RR^d}u(t,x)\psi(x)\dx\right|\le \frac{K_2}{t^{d\vartheta(m-1)}}\,.
\]
Integrating the above differential inequality we obtain for any $s,t\ge 0$:
\begin{equation}\label{est.loc.nonloc}
\int_{\RR^d}u(t,x)\psi(x)\dx\le \int_{\RR^d}u(s,x)\psi(x)\dx
+ 2s\vartheta\,K_2\,\left|t^{2s\vartheta}-s^{2s\vartheta}\right|\quad\mbox{ for all $\psi\in C_c^\infty(\RR^d)$\,.}
\end{equation}
\noindent$\bullet~$\textsc{Step 4. } We still have to pass from test functions $\psi\in C_c^\infty(\RR^d)$ to $\psi\in C_c^0(\RR^d)$  in formula \eqref{eq.trace2}, but this is easy by approximation (mollification).\qed

\medskip

We notice that the estimates \eqref{est.loc.nonloc} are only pseudo-local estimates: the global information about $u(T)$, namely the bound $\|u(T)\|_{\LL^1(\RR^d)}$ is contained in the constant  $K_2$.

\medskip


\section{Appendix I. Definitions, complements and computations}\label{sec.app}

\subsection{Definition of the fractional Laplacian.}\label{ssec.app1}

According to  Stein,   \cite{Stein70}, chapter V, the definition of the nonlocal operator
$(-\Delta)^{\sigma/2}$, known as the Laplacian of order $\sigma$, is done by means of Fourier series
\begin{equation}
  ((-\Delta)^{\sigma/2}f)^{\widehat{}}(x)=(2\pi|x|)^{\sigma} \hat f (x)\,,
\end{equation}
and can be used for positive and negative values of $\sigma$.
If $0<\sigma<2$, we can also use the representation by means of an
hypersingular kernel,
\begin{equation}\label{formula.slapl}
(-\Delta)^{\sigma/2}  g(x)= c_{d,\sigma }\mbox{
P.V.}\int_{\mathbb{R}^d} \frac{g(x)-g(z)}{|x-z|^{d+{\sigma} }}\,dz,
\end{equation}
where $c_{d,\sigma
}=\frac{2^{\sigma-1}\sigma\Gamma((d+\sigma)/2)}{\pi^{d/2}\Gamma(1-\sigma/2)}$
is a normalization constant. Another classical way of defining the fractional powers of a
linear self-adjoint nonnegative operator, in terms of the associated
semigroup, which in our case reads
\begin{equation}
\displaystyle(-\Delta)^{\sigma/2}
g(x)=\frac1{\Gamma(-\frac{\sigma}2)}\int_0^\infty
\left(e^{t\Delta}g(x)-g(x)\right)\frac{dt}{t^{1+\frac{\sigma}2}}.
\label{laplace}\end{equation}
In this paper we consistently put $\sigma=2s$, $0<s<1$ (sometimes, also $s=1$).

\subsection{Definition of weak and very weak solutions}

We recall here the definitions of weak and strong solutions taken from \cite{DPQRV2}. We finally introduce the definition of very weak solutions.

\begin{defn}\label{def:weak.solution.nonlocal} A function $u$
is a {\sl weak}  solution to Equation \eqref{FDE.eq}
if:
\begin{itemize}
\item $u\in C((0,\infty): L^1(\RR^d))$, $|u|^{m-1}u \in L^2_{\rm
loc}((0,\infty):\dot{H}^{s}(\RR^d))$;
\item  The identity
\begin{equation}
\displaystyle \int_0^\infty\int_{\RR^d}u\dfrac{\partial
\varphi}{\partial
t}\,\dx\dt-\int_0^\infty\int_{\RR^d}(-\Delta)^{s/2}(|u|^{m-1}u)(-\Delta)^{s/2}\varphi\,\dx\dt=0.
\end{equation}
holds for every $\varphi\in C_0^1(\RR^d\times(0,\infty))$;
\item  A {\sl weak} solution to Problem \eqref{FDE.eq}--\eqref{FDE.id} is a weak solution to  Equation \eqref{FDE.eq} such that moreover $u\in C([0,\infty): L^1(\RR^d))$ and $u(0,\cdot)=u_0\in \LL^1(\RR^d)$.
\end{itemize}
\end{defn}
Note that in \cite{DPQRV2} these weak solutions are given  the more precise name {\sl weak $L^1$-energy solutions}. We recall that the fractional Sobolev space $\dot{H}^{s}(\RR^d)$ is defined as the
completion of $C_0^\infty(\RR^d)$ with the norm
$$
  \|\psi\|_{\dot{H}^{s}}=\left(\int_{\RR^d}
|\xi|^\sigma|\hat{\psi}|^2\,d\xi\right)^{1/2}
=\|(-\Delta)^{s/2}\psi\|_{2}.
$$

\begin{defn} We say that a weak solution $u$  to
Problem \eqref{FDE.eq}--\eqref{FDE.id} is a strong solution if  moreover $\partial_tu\in
L^\infty((\tau,\infty):L^1(\RR^d))$, for every $\tau>0$.
\label{def:strong.solution}\end{defn}

\begin{defn}\label{def:weak.solution.nonlocal.2} A function $u$
is a {\it very weak} solution to Equation \eqref{FDE.eq} if:
\begin{itemize}
\item $u\in C((0,\infty): L^1_{\rm loc}(\RR^d))$, $|u|^{m-1}u \in L^1_{\rm loc}\left((0,\infty):L^1\left(\RR^d, (1+|x|)^{-(d+2s)}\dx\right)\right)$;
\item  The identity
\begin{equation}
\displaystyle \int_0^\infty\int_{\RR^d}u\dfrac{\partial
\varphi}{\partial
t}\,\dx\dt-\int_0^\infty\int_{\RR^d} |u|^{m-1}u\,(-\Delta)^{s}\varphi\,\dx
\dt=0.
\end{equation}
holds for every $\varphi\in C_c^\infty([0,T]\times\RR^d)$\,;
\item  A {\sl very weak} solution to Problem \eqref{FDE.eq}--\eqref{FDE.id} is very weak solution to  Equation \eqref{FDE.eq} such that moreover $u\in C([0,\infty): L^1_{\rm loc}(\RR^d))$ and $u(0,\cdot)=u_0\in \LL^1_{\rm loc}(\RR^d)$.
\end{itemize}
\end{defn}

\normalcolor

\subsection{Some functional inequalities related to the fractional Laplacian}
We recall here some useful functional inequalities which have been used throughout the paper.
\begin{lem}[Stroock-Varopoulos' inequality] Let $0<s<1$\,, $q>1$. Then
\begin{equation}\label{StrVar.ineq}
\int_{\RR^d}|v|^{q-2}v\,(-\Delta)^s v\,\dx\ge \frac{4(q-1)}{q^2}\int_{\RR^d}\left|(-\Delta)^{\frac{s}{2}}|v|^{\frac{q}{2}}\right|^s\dx\,,
\end{equation}
for all $v\in \LL^q(\RR^d)$ such that $(-\Delta)^s v\in \LL^q(\RR^d)$.
\end{lem}
\noindent\textbf{Remark. }We have used the above Stroock-Varopoulos inequality, applied to $0\le v=u^m$ and $q=(p+m-1)/m>1$\,, whenever $p>1$,  which is
\begin{equation}\label{StrVar.ineq.um}
\int_{\RR^d}|u|^{p-2}u\,(-\Delta)^s (|u|^{m-1}u)\,\dx\ge \frac{4m(p-1)}{(p+m-1)^2}\int_{\RR^d}\left|(-\Delta)^{\frac{s}{2}}|u|^{\frac{p+m-1}{2}}\right|^s\dx\,.
\end{equation}

\begin{thm}[Sobolev Inequality] Let $0<s\le 1$ and $2s<d$. Then
\begin{equation}\label{Sob.Fract}
\|f\|_{\frac{2d}{d-2s}}\le \mathcal{S}_s \left\|(-\Delta)^{s/2}f\right\|_2
\end{equation}
where the best constant is given by
\begin{equation}\label{best.const.Sob.Fract}
\mathcal{S}_s^2 := 2^{-2s}\,\pi^{-s}\frac{\Gamma\left(\frac{d-2s}{2}\right)}{\Gamma\left(\frac{d+2s}{2}\right)}
    \left[\frac{\Gamma(d)}{\Gamma(d/2)}\right]^{\frac{2s}{d}}
    =\frac{\Gamma\left(\frac{d-2s}{2}\right)}{\Gamma\left(\frac{d+2s}{2}\right)}|\mathbb{S}_d|^{-\frac{2s}{d}}
\end{equation}
and is attained on the family of functions
\[
F(x):=a\left[b^2+(x-x_0)^2\right]^{-\frac{d-2s}{2}}\,,\qquad\mbox{with $x,x_0\in\RR^d$ and $a\in \RR$\,, $b>0$\,.}
\]
\end{thm}

\subsection{Proof of Lemma \ref{Lem.phi}}\label{sec.A1}

\medskip
\noindent {\sl Proof.~}The proof is divided into several steps.

\noindent$\bullet~$\textsc{Step 1. }\textit{The integral is convergent. }First we have to prove that
\begin{equation*}
c_{d,s}^{-1}\left|(-\Delta)^s\varphi (x)\right|=\left|\int_{\RR^d}\frac{\varphi (x)-\varphi (y)}{|x-y|^{d+2s}}\dy\right|<\infty\qquad\mbox{for any }x\in \RR^d
\end{equation*}
to this end we fix $x\in \RR^d$ and we split the integral in two parts:
\[
\int_{\RR^d}\frac{\varphi (x)-\varphi (y)}{|x-y|^{d+2s}}\dy
=\int_{|x-y|>\delta}\frac{\varphi (x)-\varphi (y)}{|x-y|^{d+2s}}\dy
+\int_{|x-y|\le \delta}\frac{\varphi (x)-\varphi (y)}{|x-y|^{d+2s}}\dy= I + II
\]
where $\delta>0$ is taken so small that the following Taylor expansion around $x\in\RR^d$ holds true
\[
\varphi (y)= \varphi (x)+\nabla\varphi (x)\cdot(y-x)+(y-x)^t \,{\rm D^2}\varphi (\overline{x})\,(y-x)
\]
for some $\overline{x}\in B_1(x)$\,. Therefore we have
\[
\begin{split}
I &= \left|\int_{|x-y|\le \delta}\frac{\varphi (x)-\varphi (y)}{|x-y|^{d+2s}}\dy\right|\\
&= \left|\int_{|x-y|\le \delta}\frac{\nabla\varphi (x)\cdot(y-x)}{|x-y|^{d+2s}}\dy +
\int_{|x-y|\le \delta}\frac{(y-x)^t \,{\rm D^2}\varphi (\overline{x})\,(y-x)}{|x-y|^{d+2s}}\dy\right|\\
&\le_{(a)} \sup_{1\le i,j\le d}\|\partial_{ij}\varphi \|_{\LL^\infty(\RR^d)} \left|\int_{|x-y|\le \delta}\frac{1}{|x-y|^{d-(2-2s)}}\dy\right|\\
&\le \sup_{1\le i,j\le d}\|\partial_{ij}\varphi \|_{\LL^\infty(\RR^d)} \int_0^{\delta}\frac{\rd r}{r^{1-2(1-s)}}
=_{(b)}K\frac{\delta^{2(1-s)}}{(2(1-s))}
\end{split}
\]
where in $(a)$ we have used that
\[
P.V.\int_{|x-y|\le \delta}\frac{\nabla\varphi (x)\cdot(y-x)}{|x-y|^{d+2s}}\dy=0
\]
for symmetry reasons. In $(b)$ we used the fact that $|\partial_{ij}\varphi (z)|\le K$ for some positive constant $K$ that depends only on $\alpha$\,. On the other hand, the outer integral is easily seen to be finite, indeed
\[
II= \left|\int_{|x-y|>\delta}\frac{\varphi (x)-\varphi (y)}{|x-y|^{d+2s}}\dy\right|
\le 2\|\varphi \|_{\LL^\infty(\RR^d)}\left|\int_{|x-y|>\delta}\frac{1}{|x-y|^{d+2s}}\dy\right|
\le 2\omega_d \int_\delta^{\infty}\frac{\rd r}{r^{1+2s}}=\frac{\omega_d}{s\delta^{2s}}\,.
\]
The above estimates for $I$ and $II$ do not depend on $x\in\RR^d$, hence $\left|(-\Delta)^s\varphi (x)\right|$ is finite for all $x\in \RR^d$\,.
\begin{figure}[ht]
\centering
\ifpdf
    \includegraphics[height=6cm, width=6cm]{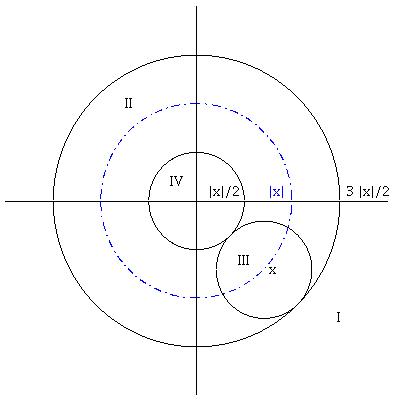}
\else
    \includegraphics[height=6cm, width=6cm]{s-testfunct.eps}
\fi
 \caption{\noindent\textit{The $4$ regions in which we split the integral}}
 \label{fig.3}
\end{figure}

\noindent$\bullet~$\textsc{Step 2. }\textit{Better estimates for $|x|$ large. } We are going to use the hypothesis that $\varphi$ is radially symmetric and decreasing for $|x|\ge 1$ and that $\varphi(x)\le |x|^{-\alpha}$\,, $|D^2\varphi(x)| \le c_0 |x|^{-\alpha-2}$\,, for some positive constant $\alpha$ and for $|x|$ large enough. We are interested in the behaviour of $|(-\Delta)^s\varphi (x)|$ for large values of $x$, therefore we fix $x\in \RR^d$ with $|x|$ sufficiently large.  We have to estimate
\begin{equation*}
c_{d,s}^{-1}\left|(-\Delta)^s\varphi (x)\right|=\left|\int_{\RR^d}\frac{\varphi (x)-\varphi (y)}{|x-y|^{d+2s}}\dy\right|\,,
\end{equation*}
to this end we split the integral into four parts, see Figure \ref{fig.3},
\begin{equation}\label{splitting.4}\begin{split}
\int_{\RR^d}\frac{\varphi (x)-\varphi (y)}{|x-y|^{d+2s}}\dy
=\int_{|y|> 3|x|/2}\frac{\varphi (x)-\varphi (y)}{|x-y|^{d+2s}}\dy
+\int_{\left\{|x|\le 2|y|\le 3|x|\right\}\setminus B_{|x|/2}(x)}\frac{\varphi (x)-\varphi (y)}{|x-y|^{d+2s}}\dy\\
+\int_{ B_{|x|/2}(x)}\frac{\varphi (x)-\varphi (y)}{|x-y|^{d+2s}}\dy
+\int_{|y|<|x|/2}\frac{\varphi (x)-\varphi (y)}{|x-y|^{d+2s}}\dy
=I+II+III+IV
\end{split}
\end{equation}
We estimate the four integrals separately, keeping in mind that we are assuming $\varphi\ge 0$ in this latter case. The first integral can be estimated as follows
\[
I=\int_{|y|> 3|x|/2}\frac{\varphi (x)-\varphi (y)}{|x-y|^{d+2s}}\dy
\le \omega_d \varphi (x) \int_{\frac{3|x|}{2}}^{\infty}\frac{\rd r}{r^{1+2s}}=\frac{k_1}{|x|^{\alpha+2s}}
\]
since $\varphi (y)\le \varphi (x)$ when $|y|> 3|x|/2$, therefore $|\varphi (x)-\varphi (y)|\le \varphi (x)$\,, and we remark that the constant $k_1$ depends only on $\alpha, s, d$, since $\varphi(x)\le |x|^{-\alpha}$ and $|x|$ is large enough. The second integral gives
\[
II\le \left|\int_{\left\{|x|\le 2|y|\le 3|x|\right\}\setminus B_{|x|/2}(x)}\frac{\varphi (x)-\varphi (y)}{|x-y|^{d+2s}}\dy\right|
\le \frac{\varphi (x/2) }{\big(|x|/2\big)^{d+2s}}\int_{\frac{|x|}{2}}^{\frac{3|x|}{2}}r^{d-1}\rd r\le  \frac{k_2}{|x|^{\alpha+2s}}
\]
since $\varphi (y)\le \varphi (x/2)$ when $|y|> |x|/2$, therefore $|\varphi (x)-\varphi (y)|\le \varphi (x/2)$\,, and we remark that the constant $k_2$ depends only on $\alpha, s, d$, since $\varphi(x/2)\le |x|^{-\alpha}$ and $|x|$ is large enough.\\
We can estimate the third integral as follows:
\[
\begin{split}
III &= \left|\int_{|x-y|\le \frac{|x|}{2}}\frac{\varphi (x)-\varphi (y)}{|x-y|^{d+2s}}\dy\right|\\
&= \left|\int_{|x-y|\le \frac{|x|}{2}}\frac{\nabla\varphi (x)\cdot(y-x)}{|x-y|^{d+2s}}\dy +
\int_{|x-y|\le \frac{|x|}{2}}\frac{(y-x)^t \,{\rm Hess}\varphi (\overline{x})\,(y-x)}{|x-y|^{d+2s}}\dy\right|\\
&\le_{(a)} \sup_{1\le i,j\le d}\|\partial_{ij}\varphi \|_{\LL^\infty(B_{\frac{|x|}{2}}(x))} \left|\int_{|x-y|\le \frac{|x|}{2}}\frac{1}{|x-y|^{d-(2-2s)}}\dy\right|\\
&\le \frac{k'_3}{|x|^{\alpha+2}} \int_0^{\frac{|x|}{2}}\frac{\rd r}{r^{1-2(1-s)}}
=\frac{k'_3}{|x|^{\alpha+2}}\left(\frac{|x|}{2}\right)^{2-2s}=\frac{k_3}{|x|^{\alpha+2s}}
\end{split}
\]
where in $(a)$ we have used that
\[
P.V.\int_{|x-y|\le \delta}\frac{\nabla\varphi (x)\cdot(y-x)}{|x-y|^{d+2s}}\dy=0
\]
for symmetry reasons as in Step 2. In $(b)$ we used the fact that $|z-x|<|x|/2$ implies $|x|/2<|z|<3|x|/2$, therefore $|\partial^2_{ij}\varphi (z)|\le c_0/|z|^{\alpha+2}\le 2^{\alpha+2}c_0/|x|^{\alpha+2}$ for all $z\in B_{|x|/2}(x)$, recalling that $|x|$ is always taken large enough. The constants $k'_3$ and $k_3$ depend only on $\alpha, s,d$.\\
It only remains to estimate the fourth integral:
\[
IV \le \int_{|y|<|x|/2}\frac{\big|\varphi (x)-\varphi (y)\big|}{|x-y|^{d+2s}}\dy
\le \frac{2^{d+2s}}{|x|^{d+2s}}\int_{|y|<|x|/2}\varphi (y)\dy.
\]
since we observe that $|y|<|x|/2$ implies $\varphi (x)\le \varphi (2y)\le \varphi (y)$ which gives $|\varphi (x)-\varphi (y)|\le \varphi (y)$, moreover we also have that $|y|<|x|/2$ implies that $|y-x|> |x|/2$. The term represents the long-range influence of the inner core of the function at large distances and will make for different conclusions of the lemma depending on the case. Indeed, we have the following estimates for $|x|$ large enough:\begin{itemize}
\item If $\alpha>d$ the last integral is finite and we get $IV\le k_4/|x|^{d+2s}$.
\item If $\alpha< d$ the last integral grows like $|x|^{d-\alpha}$ and we get $IV\le k_5/|x|^{\alpha+2s}$.
\item Finally when $\alpha=d$ we get $IV\le k_6\log|x|/|x|^{d+2s}$\,.
\end{itemize}
We finally remark that the constants $k_4,k_5,k_6$ depend only on $\alpha,s,d$\,.

\noindent$\bullet~$\textsc{Step 3. }\textit{Positivity estimates for $|x|$ large. }In the case when $\alpha>d$ we need to prove that if $\varphi\ge0$ then we have that $|(-\Delta)^s\varphi(x)|\ge c_4 |x|^{-(d+2s)}$  for all $|x|\ge |x_0|>>1$\,. We split the integral into four parts, as in Step 2, equation \eqref{splitting.4}, see Figure \ref{fig.3},
\[
c_{d,s}^{-1}(-\Delta)^s\varphi (x)=\int_{\RR^d}\frac{\varphi (x)-\varphi (y)}{|x-y|^{d+2s}}\dy
=I+II+III+IV
\]
We have proven in Step 2 that $|I|+|II|+|III|\le (k_1+k_2+k_3)\,/|x|^{\alpha+2s}$ and we recall that the constant $k_i$ depend only on $\alpha, s, d$. We just have to obtain better estimates for the last term, to this end we further split the integral in two parts:
\[\begin{split}
IV = \int_{|y|<|x|/2}\frac{\varphi (x)-\varphi (y)}{|x-y|^{d+2s}}\dy
=\int_{|y|<|x|/2}\frac{\varphi (x)}{|x-y|^{d+2s}}\dy -\int_{|y|<|x|/2}\frac{\varphi (y)}{|x-y|^{d+2s}}\dy = IV_a- IV_b
\end{split}
\]
Let us calculate
\[
0\le IV_a=\int_{|y|<|x|/2}\frac{\varphi (x)}{|x-y|^{d+2s}}\dy
\le \frac{\varphi(x)}{(|x|/2)^{d+2s}} \int_{|y|<|x|/2}\dy
\le \frac{k_4}{|x|^{\alpha+2s}}
\]
since $|x-y|\ge |x|/2$ when $|y|\le |x|/2$ and $\varphi(x)\le |x|^{-\alpha}$. We remark that the constant $k_4$ depends only on $\alpha, s, d$. On the other hand, $IV_b\ge 0$ and
\[
0\le IV_b=\int_{|y|<|x|/2}\frac{\varphi (y)}{|x-y|^{d+2s}}\dy
\le \frac{1}{(|x|/2)^{d+2s}} \int_{|y|<|x|/2}\varphi(y)\dy
\le \frac{\|\varphi\|_{\LL^1(\RR^d)}}{|x|^{d+2s}}
\]
Summing up, we have obtained that
\[
\begin{split}
-(-\Delta)^s\varphi (x)&=-c_{d,s}\int_{\RR^d}\frac{\varphi (x)-\varphi (y)}{|x-y|^{d+2s}}\dy
\ge c_{d,s}\left[IV_b -\big(\,|I|+|II|+|III|+|IV_a|\,\big)\right]\\
&\ge c_{d,s}\left[\frac{\|\varphi\|_{\LL^1(\RR^d)}}{|x|^{d+2s}}-\frac{k_5}{|x|^{\alpha+2s}}\right]
=\left[\|\varphi\|_{\LL^1(\RR^d)}-\frac{k_5}{|x|^{\alpha-d}}\right]\frac{c_{d,s}}{|x|^{d+2s}}\ge \frac{c_4}{|x|^{d+2s}}
\end{split}
\]
since $|I|+|II|+|III|+|IV_a|\le (k_1+k_2+k_3+k_4)\,|x|^{\alpha+2s}=k_5\,|x|^{\alpha+2s}$ and $c_4>0$ since $\alpha>d$, if we choose $|x|$ sufficiently large, namely $|x|^{\alpha-d}\ge k_5/\|\varphi\|_{\LL^1(\RR^d)}$\,.\qed

\normalcolor

\subsection{Optimization Lemma}\label{app.opt}

We state and prove here a simple technical lemma  that has been used in the proof of Theorem \ref{thm.lower}\,.

\begin{lem}\label{Lem.Opt} Let $0<m,s<1$\,, $2s>d(1-m)$\,, $\vartheta=1/[2s-d(1-m)]>0$ and $B,C,t>0$. Define
\[
F(t,R):=\frac{A(t)}{R^{d(1-m)}}-\frac{B\,t}{R^{2s}}\,,\qquad\mbox{with}\qquad A(t):= M-\frac{C}{t^{d(1-m)\vartheta}}\,.
\]
Then there exists
\begin{equation}\label{t.min}
t_*:=2s\vartheta\left(\frac{C}{M}\right)^{\frac{1}{d(1-m)\vartheta}}>0
\end{equation}
and
\begin{equation}\label{R.max}
\overline{R}(t)=\left(\frac{2sBt}{d(1-m)A(t)}\right)^{\vartheta}
\ge\overline{R}(t_*)=\left[\frac{2s}{d(1-m)}\frac{(2s\vartheta)^{2s\vartheta}}{(2s\vartheta)^{d(1-m)}-1}\right]^{\vartheta}\frac{B^\vartheta C^{\frac{1}{d(1-m)}}}{M^{\frac{2s\vartheta}{d(1-m)}}}
>0
\end{equation}
so that for all $t\ge t_*$ we have
\[
F(\overline{R}(t),t)=\left[\left(\frac{2s}{d(1-m)}\right)^{\frac{1}{\vartheta}}-1\right]\left[\frac{d(1-m)}{2s}\right]^{2s\vartheta}
\frac{A(t)^{2s\vartheta}}{(Bt)^{d(1-m)\vartheta}}>0\,.
\]
\end{lem}
\noindent {\sl Proof.~}First we observe that $A(t)$ is monotone increasing in $t>0$, and that $A(t_*/2s\vartheta)=0$, where $t_*$ has the expression given by \eqref{t.min}, so that $A(t)>A(t_*)>A(t/2s\vartheta)=0$ since $2s\vartheta>1$, and
\[
A(t_*)=\frac{(2s\vartheta)^{d(1-m)\vartheta}-1}{(2s\vartheta)^{d(1-m)\vartheta}}M>0
\]
Moreover, it easy to check that $t_*$ is also the value for which $A(t_*)-t_*A'(t_*)=0$. Next we fix a time $t\ge t_*$ and we find the maximum with respect to $R$ of the function $F(t,R)$:
\[
\partial_R F(R,t)=-\frac{d(1-m)A(t)}{R^{d(1-m)\vartheta+1}}+\frac{2sB\,t}{R^{2s+1}}\,,
\]
and $\partial_R F(\overline{R}(t),t)=0$, so that the maximum is attained at $\overline{R}(t)$ whose expression it is easily checked to be \eqref{R.max}. It only remains to prove that $\overline{R}(t)\ge\overline{R}(t_*)>0$, to this end we observe that
\[
\partial_t \overline{R}(t)=\vartheta\left[\frac{2sB}{d(1-m)}\right]^\vartheta\left[\frac{t}{A(t)}\right]^{\vartheta-1}\frac{A(t)-tA'(t)}{A(t)^2}
\]
and it is clear now that the minimum is attained at $t_*$\,, since $\partial_t \overline{R}(t_*)=0$\,, because we already know that  $A(t_*)-t_*A'(t_*)=0$.\qed

\subsection{Reminder about measure theory}\label{appendix.meas}

We recall here some basic facts on measure theory for convenience of the reader. We refer the interested reader to the books \cite{EvansMeas, RudinRealCplx}.

\begin{defn}
A measure $\mu$ is \textsl{regular} if
\[
\forall A\subseteq \RR^d \; \exists B\;\mu\mbox{-measurable such that }A\subseteq B\; \mbox{and}\;\mu(A)=\mu(B)\,.
\]
A measure $\mu$ is \textsl{Borel} if every Borel set $\mathcal{B}(\RR^d)$ is $\mu$-measurable. A measure $\mu$ is \textsl{Borel regular }if
\[
\forall A\subseteq \RR^d \; \exists B\in\mathcal{B}(\RR^d)\mbox{  such that }A\subseteq B\; \mbox{and}\;\mu(A)=\mu(B)\,.
\]
A measure $\mu$ is \textsl{Radon} if is Borel regular and $\mu(K)<+\infty$ for any compact set $K\subset\RR^d$.\\
A sequence of measures $\mu_n$ \textsl{converges weakly (star)} to the measure $\mu$, $\mu_n\rightharpoonup \mu$ as $n\to\infty$ if
\[
\lim_{n\to\infty}\int_{\RR^d}\varphi\,\rd\mu_n=\int_{\RR^d}\varphi\,\rd\mu\qquad\mbox{for all }\varphi\in C^0_c(\RR^d)\,.
\]
\end{defn}
\begin{thm}[Weak compactness for measures]\label{thm.cpt.radon}
Let $\{\mu_n\}$ be a sequence of Radon measures on $\RR^d$ satisfying
\[
\sup_{n}\mu_n(K)<\infty\qquad\mbox{for any compact set $K\subset\RR^d$.}
\]
Then there exists a subsequence $\mu_{n_k}$ and a Radon measure $\mu$ such that $\mu_{n_k}\rightharpoonup \mu$ as $k\to\infty$\,.
\end{thm}


\begin{thm}[Riesz Representation Theorem]\label{thm.riesz.radon}
Assume $L:C_c^{\infty}(\RR^d)\to\RR$ is linear and nonnegative, so that
\[
L\varphi\ge 0\qquad\mbox{for all }0\le \varphi\in C_c^{\infty}
\]
Then there is a unique Radon measure $\mu$ on $\RR^d$ such that
\[
L\varphi=\int_{\RR^d}\varphi\,\rd\mu \qquad\mbox{for all } \varphi\in C_c^{\infty}
\]
\end{thm}


\color{darkblue}
\section{Appendix II. Applied literature and motivation}\label{App.Motiv}

We gather here some updated information on the occurrence of the nonlinear fractional diffusion equation we propose and related models in the physical or probabilistic literature.

 $\bullet$  A great variety of diffusive problems in nature, namely those referred to as normal diffusion, are satisfactorily described by the classical Heat Equation or Fokker-Planck linear equation. However, anomalous diffusion is nowadays intensively studied, both theoretically and experimentally since it conveniently explains a number of phenomena in several areas of physics, biology, ecology, geophysics, and many others, which can be  briefly summarized as having non-Brownian scaling. A large variety of phenomena in physics and finance are modeled by linear anomalous diffusion equations, see e.g. \cite{A2009, C2004, W2001}. Fractional kinetic equations of the diffusion, diffusion-advection, and Fokker-Planck type represent a useful approach for the description of transport dynamics in complex systems which are governed by anomalous diffusion. These fractional equations are usually derived asymptotically from basic random walk models, cf.  \cite{MK2000}.

 $\bullet$  Anomalous diffusion often takes a nonlinear form. To be more specific, there exist many phenomena  in nature where, as time goes on, a crossover is observed between different diffusion regimes. Tsallis et al. \cite{BGT2000, LMT2003}  discuss the following cases: (i) a mixture of the porous medium equation, which is connected with non-extensive statistical mechanics, with the normal diffusion equation; (ii) a mixture of the fractional time derivative and normal diffusion equations; (iii) a mixture of the fractional space derivative, which is related with L\'evy flights, and normal diffusion equations. In all three cases a crossover is obtained between anomalous and normal diffusions. This leads to models of nonlinear diffusion of porous medium or fast diffusion types with standard or  fractional Laplace operators, cf.  eqn. (4) of \cite{BGT2000}.

$\bullet$  There have been many studies of hydrodynamic limits of interacting particle systems with long-range dynamics, which lead to fractional diffusion equations of our type, mainly linear like in \cite{JKOlla}, \cite{Jara0}, but also nonlinear in the recent literature, cf. the works \cite{Jara1}, \cite{Jara2}. Thus, in the last reference, Jara and co-authors study the non-equilibrium functional central limit theorem for the position of a tagged particle in a mean-zero one-dimensional zero-range process. The asymptotic
behavior of the particle is described by a stochastic differential equation
governed by the solution of the following nonlinear hydrodynamic (PDE) equation, $\partial_t \rho = \sigma^2 \partial^2_x \Phi(\rho)$. When $\Phi$ is a power we recover equation \eqref{FDE.eq}.

$\bullet$ Equations like the last one (in several space dimensions) occur in boundary heat control, as already mentioned by Athanasopoulos and Caffarelli \cite{AC2010}\,, where they refer to the model  formulated in the book by Duvaut and Lions \cite{DL1972}, and use the so-called Caffarelli-Silvestre extension \cite{CS2007}.

$\bullet$ The combination of diffusion with convection is an important research topic with abundant literature. The use of nonlinear fractional diffusion in that setting has been studied by Cifani and Jakobsen \cite{CJ}, where references to other models is given.

$\bullet$ A different version of the nonlinear fractional diffusion equation of porous medium type takes the form $u_t=\nabla (u\, \nabla {\cal K} u)$, where $K$ is the Riesz operator that expresses the inverse to the fractional Laplacian, ${\cal K} u=(-\Delta)^{-s}u$. This has been studied by Caffarelli and Vazquez in \cite{CV2010, CV2011} and Biler, Karch and Monneau in \cite{BKM2010}, where the equation is derived in the framework of the theory of dislocations.
This model has some strikingly different properties, like lack of strict positivity and occurrence of free boundaries. See the survey \cite{VazAbel}.
\normalcolor

\newpage

\noindent {\textbf{\large \sc Acknowledgment.} Both authors funded have been partially funded by Project MTM2011-24696 (Spain). We are grateful to A. de Pablo, F. Quir\'os and A. Rodriguez for discussions of the topic, and to the referee for enlightening suggestions.

\vskip .3cm



\begin{thebibliography}{00}
\small
%
%
%
%
\color{darkblue}
\bibitem{A2009} {\sc D. Applebaum, } ``L\'evy Processes and Stochastic Calculus'', second edition, Cambridge Stud. Adv. Math., vol. 116, Cambridge University Press, Cambridge, 2009.
\normalcolor


\bibitem{ArCaff} {\sc  D.~G. Aronson, L.~A. Caffarelli}, The initial
trace of a solution of the porous medium equation, {\em Trans. Amer.
Math. Soc.} \textbf{280} (1983), 351--366.


\color{darkblue}
\bibitem{AC2010} {\sc I. Athanasopoulos, L. A. Caffarelli, }Continuity of the temperature in boundary heat control problems. \textit{Adv. Math. } \textbf{224} (2010), no. 1, 293–315.

\bibitem{BPSV2013} {\sc B. Barrios, I. Peral, F. Soria, E. Valdinoci, }A Widder's type Theorem for the heat equation with nonlocal diffusion, \textit{Preprint }(2013), \texttt{http://arxiv.org/abs/1302.1786}
\normalcolor

\bibitem{BCr-cont} {\sc  P. B\'enilan, M.~G. Crandall, }The continuous dependence on $\varphi $  of solutions of \ $u_t - \Delta \varphi (u) = 0$, {\em  Indiana Univ. Math. J.}  {\bf  30} (1981),   161--177

\bibitem {BCr} {\sc  P. B\'enilan, M.~G. Crandall,} Regularizing effects of homogeneous evolution equations. {\em Contributions to Analysis and Geometry} (suppl.  to Amer. Jour. Math.), Johns Hopkins Univ. Press, Baltimore, Md., 1981. Pp.
23-39.

\color{darkblue}
\bibitem{BKM2010} {\sc Biler, P., Karch, G.,  Monneau, R.}, {Nonlinear diffusion of dislocation density and self-similar
solutions}. \textit{Comm. Math. Phys.} \textbf{294} (2010), no.~1, 145--168.

\bibitem{BGT2000} {\sc M. Bologna,  P. Grigolini, C. Tsallis, }Anomalous diffusion associated with nonlinear fractional derivative Fokker-Planck-like equation: Exact time-dependent solutions, \textit{Physical Review E }\textbf{62}, (2000).
\normalcolor

\bibitem{BV} {\sc M.  Bonforte, J.~L. V\'azquez, }Global positivity estimates and {H}arnack inequalities for the fast diffusion equation, \textit{J. Funct. Anal.} \textbf{240}, no.2, (2006) 399--428.

\bibitem{BV3} {\sc M.  Bonforte, J.~L. V\'azquez, }Reverse Smoothing Effects, Fine Asymptotics and Harnack Inequalities for Fast Diffusion Equations\,, \textit{Bound. Value Prob.} Volume 2007, Article ID 21425, 31 pages, doi:10.1155/2007/21425

\bibitem{BV-ADV} {\sc M.  Bonforte, J.~L. V\'azquez, }{\rm Positivity, local smoothing, and Harnack inequalities for very fast diffusion equations}, \textit{Advances in Math.} \bf 223 \rm (2010), 529--578.


\color{darkblue}


\bibitem{CS2007} {\sc L.A. Caffarelli, L. Silvestre, }An extension problem related to the fractional Laplacian, \textit{Comm. Partial Differential Equations,} \textbf{32} (2007) 1245–1260.

\bibitem{CV2011} {\sc Caffarelli, L.~A., V\'azquez, J.~L.}, {Nonlinear porous medium flow with fractional potential pressure}.
\textit{Arch. Rational Mech. Anal.} \textbf{202} (2011) 537–-565.

\bibitem{CV2010}  {\sc Caffarelli, L.~A.,  V\'azquez, J.~L.}, {Asymptotic behaviour of a porous medium equation with fractional
diffusion}. \textit{Discrete Contin. Dyn. Syst.} \textbf{29} (2011), no.~4, 1393--1404.
\normalcolor


\bibitem {ChVaz02} {\sc E. Chasseigne, J.~L. V\'azquez},  Extended theory of fast diffusion equations in optimal classes of data. Radiation from singularities, \textit{Arch. Rat. Mech. Anal.}  {\bf 164} (2002), 133--187.

\color{darkblue}
\bibitem{CJ}  {\sc Cifani, S.; Jakobsen, E.~R.}, \textit{Entropy solution theory for fractional degenerate convection-diffusion equations}. \textrm{Ann. Inst. H. Poincar\'{e} Anal. Non Lin\'{e}aire} \textbf{28} (2011), no. 3, 413–-441.

\bibitem{C2004} {\sc R. Cont, P. Tankov, } {\em ``Financial Modelling with Jump Processes'', } Chapman \& Hall/CRC, 2004.
\normalcolor

\bibitem{DK2} {\sc B. Dahlberg, C.~E. Kenig}, Nonnegative solutions to fast diffusions, \textit{Rev. Mat. Iberoamericana} \textbf{4} (1988), n. 1, 11-–29.

\color{darkblue}


\bibitem{DBbook} {\sc E. DiBenedetto}. {\em ``Degenerate parabolic equations'',} Universitext. Springer-Verlag, New York, 1993.

\bibitem{DGVbook} E. DiBenedetto, U. Gianazza, V. Vespri, \textrm{``Harnack'\,s inequality for degenerate and singular parabolic equations''\,, } Springer Monographs in Mathematics, Springer 2011.

\bibitem{DL1972} {\sc G. Duvaut, J.-L. Lions, }{\em ``Les In{\'e}quations en Mechanique et en Physique''}, Travaux et Recherches Mathématiques, No. 21. Dunod, Paris, 1972. xx+387 pp.
\normalcolor

\bibitem{EvansMeas} {\sc L.~C. Evans, R.~F. Gariepy, } Measure theory and fine properties of functions. \textit{Studies in Advanced Mathematics,} CRC Press, Boca Raton, FL, 1992. viii+268 pp. ISBN: 0-8493-7157-0

\bibitem{HP} {\sc M.~A. Herrero, M. Pierre, }{The Cauchy Problem for $u_t=\Delta u^m$ when $0<m<1$}, \textit{Trans. Amer. Math. Soc. }\textbf{291} n.1 (1985), 145--158.

\color{darkblue}
\bibitem{Jara0} {\sc M. D. Jara,} Nonequilibrium scaling limit for a tagged particle in the simple exclusion process with long jumps, \textit{Comm. Pure Applied Math.}, {\bf 62} (2009),  198--214.


\bibitem{Jara1} {\sc M. Jara}, Hydrodynamic Limit Of Particle Systems With Long Jumps\,, {\it Preprint, }(2009). \\ \texttt{http://arxiv.org/abs/0805.1326v2}

\bibitem{JKOlla} {\sc M. D. Jara, T. Komorowski, S.  Olla,} {Limit theorems for additive functionals of a Markov chain}. \textit{Ann. Appl. Probab.} \textbf{19} (2009), no. 6, 2270–-2300.

\bibitem{Jara2} {\sc  Jara, C. Landim,  S. Sethuraman, } Nonequilibrium fluctuations for a tagged particle in mean-zero one-dimensional zero-range processes, \textit{ Probab. Theory Relat. Fields }{\bf 145 }(2009), 565--590.

\bibitem{LMT2003} {\sc E. K. Lenzi, R. S. Mendes, C. Tsallis, }Crossover in diffusion equation: Anomalous and normal behaviors, \textit{Physical Review E }\textbf{67}, 031104 (2003).

\bibitem{MK2000} {\sc R. Metzler, J. Klafter}, The random walk's guide to anomalous diffusion: a fractional dynamics approach, {\sl  Physics Reports} {\bf 339} (2000) 1-77.
\normalcolor

\bibitem{DPQRV1} {\sc A.~de Pablo, F. Quir{\'o}s, A. Rodriguez, J. L. V´azquez, }\textrm{A fractional porous medium equation} Adv. Math. \textbf{226} (2011), no. 2, 1378-–1409.

\bibitem{DPQRV2} {\sc A.~de Pablo, F. Quir{\'o}s, A. Rodriguez, J. L. V´azquez, }  \textrm{A general fractional porous medium equation}, \textit{Comm. Pure Appl. Math. }\textbf{65} (2012), no. 9, 1242-–1284.

\bibitem{Pierre} {\sc  M. Pierre,}{ Uniqueness of the solutions of $u_t-\Delta \phi(u)=0$ with initial datum a measure}. \textit{ Nonlinear Anal. T. M. A.}  {\bf 6} (1982), pp. 175--187.

\color{darkblue}
\bibitem{PT2}{\sc F. Punzo, G. Terrone, }On the Cauchy problem for a general fractional porous medium equation with variable density, {\it Preprint, }(2013). \texttt{http://arxiv.org/abs/1302.0119}
\normalcolor




\bibitem{RudinRealCplx} {\sc W. Rudin, } {\em ``Real and complex analysis. Third edition''.} \textit{McGraw-Hill Book Co.}, New York, 1987. xiv+416 pp. ISBN: 0-07-054234-1

\color{darkblue}
\bibitem{SV2013} {\sc  D. Stan, J.~L. V{\'a}zquez, } The Fisher-KPP equation with nonlinear fractional diffusion, {\em Preprint, }(2013). \texttt{http://arxiv.org/abs/1303.6823}
\normalcolor

\bibitem{Stein70} {\sc E. Stein.} {\em  ``Singular Integrals and Differentiability Properties of Functions''}, Princeton University Press, Princeton, 1970.


\bibitem{VazLN} {\sc  J.~L. V{\'a}zquez, }{\em ``Smoothing and decay estimates for nonlinear diffusion equations''}, vol.~33 of Oxford Lecture Notes in Maths. and its Applications, Oxford Univ. Press, 2006.

\bibitem{VazBook} {\sc  J.~L. V{\'a}zquez, } {\em ``The Porous Medium Equation. Mathematical Theory''}, vol.~Oxford Mathematical Monographs, Oxford University Press, Oxford, 2007.

\bibitem{Vaz2012} {\sc  J.~L. V{\'a}zquez, } Barenblatt solutions and asymptotic behaviour for a nonlinear fractional heat equation of porous medium type, {\it To appear in J. Europ. Math. Soc., }(2013). Preprint: \texttt{http://arxiv.org/abs/1205.6332}

\color{darkblue}
\bibitem{VazAbel} {\sc V\'{a}zquez, J.~L.}, {Nonlinear diffusion with fractional Laplacian operators}. In \textit{Nonlinear Partial Differential Equations: The Abel Symposium 2010},  pp. 271--298, Abel Symposia, 7, Springer-Verlag, Berlin, 2012.

\bibitem{VV2013} {\sc  J.~L. V{\'a}zquez, B. Volzone, }Symmetrization for Linear and Nonlinear Fractional Parabolic Equations of Porous Medium Type, {\it Preprint, }(2013). \texttt{http://arxiv.org/abs/1303.2970}

\bibitem{W2001} {\sc W. Woycz\'ynski, } \textit{Lévy processes in the physical sciences, }in: Lévy Processes, Birkhäuser, Boston, 2001, pp. 241
\normalcolor


\end{thebibliography}
\end{document}